\documentclass[11pt,reqno]{amsart}
\usepackage{amssymb,verbatim}
\newcommand{\la}{\lambda}
\newcommand{\ssp}{s^{\vdom}}
\newcommand{\so}{s^{\hdom}}
\newcommand{\sd}{s^{\cell}}
\newcommand{\gggg}{\mathfrak{g}}
\newcommand{\gb}{\overline{\gggg}}
\newcommand{\Xb}{\overline{X}}
\newcommand{\Mb}{\overline{M}}
\newcommand{\se}{\diamondsuit}
\newcommand{\nn}{\varnothing}
\newcommand{\PPP}{\mathcal{P}}
\newcommand{\NN}{\mathbb{Z}_{\ge0}}
\newcommand{\La}{\Lambda}

\newcommand{\End}{\mathrm{End}}
\newcommand{\Vt}{\widetilde{V}}

\newcommand{\R}{\mathbb{R}}
\newcommand{\ac}{\alpha^\vee}
\newcommand{\Rb}{\overline{R}}
\newcommand{\rrho}{{\overline{\rho}}}

\newcommand{\Z}{\mathbb{Z}}
\newcommand{\N}{\mathbb{Z}_{\ge0}}
\newcommand{\HH}{\mathbb{H}}

\newcommand{\tB}{\widetilde{B}}
\newcommand{\tBB}{\mathbb{B}}
\newcommand{\Par}{\mathcal{P}}

\newcommand{\KK}{f}
\newcommand{\minus}{\varepsilon}
\newcommand{\sign}{\epsilon}

\newcommand{\coeff}{\Big|}

\newcommand{\inner}[2]{\left\langle #1\,,\,#2\right\rangle}

\newdimen\squaresize \squaresize=4pt
\newdimen\thickness \thickness=0.2pt

\def\square#1{\hbox{\vrule width \thickness
     \vbox to \squaresize{\hrule height \thickness\vss
        \hbox to \squaresize{\hss#1\hss}
     \vss\hrule height\thickness}
\unskip\vrule width \thickness}
\kern-\thickness}

\def\vsquare#1{\vbox{\square{$#1$}}\kern-\thickness}

\def\young#1{
\vbox{\smallskip\offinterlineskip
\halign{&\vsquare{##}\cr #1}}}

\def\thisbox#1{\kern-.09ex\fbox{#1}}
\def\downbox#1{\lower1.200em\hbox{#1}}
\newcommand{\hdom}{{\begin{picture}(8,4)\multiput(0,0)(4,0){3}{\line(0,1)
{4}}%
\multiput(0,0)(0,4){2}{\line(1,0){8}}\end{picture}}}
\newcommand{\vdom}{{\begin{picture}(4,8)\multiput(0,0)(0,4){3}{\line(1,0)
{4}}%
\multiput(0,0)(4,0){2}{\line(0,1){8}}\end{picture}}}
\newcommand{\cell}{{\begin{picture}(4,4)\multiput(0,0)(4,0){2}{\line(0,1)
{4}}%
\multiput(0,0)(0,4){2}{\line(1,0){4}}\end{picture}}}

\newcommand{\cclccl}{{\begin{picture}(8,8)\multiput(0,0)(4,0){3}{\line(0,1){8}}%
\multiput(0,0)(0,4){3}{\line(1,0){8}}\end{picture}}}

\newtheorem{theorem}{Theorem}
\newtheorem{prop}[theorem]{Proposition}
\newtheorem{conjecture}[theorem]{Conjecture}
\newtheorem{lemma}[theorem]{Lemma}
\newtheorem{cor}[theorem]{Corollary}
\theoremstyle{definition}
\newtheorem{remark}[theorem]{Remark}

\newtheorem{example}[theorem]{Example}

\numberwithin{equation}{section}

\begin{document}

\title[Deformed universal characters]{Deformed universal characters \\
for classical and affine algebras}

\author{Mark Shimozono}
\address{Department of Mathematics \\
Virginia Tech \\
Blacksburg, VA 24061-0123 USA}
\email{mshimo@math.vt.edu}
\thanks{First author is partially supported by NSF grant DMS-0100918}.
\author{Mike Zabrocki}

\begin{abstract} Creation operators are given for the
three distinguished bases of the type $BCD$ universal character
ring of Koike and Terada yielding an elegant way of treating
computations for all three types in a unified manner.
Deformed versions of these operators
create symmetric function bases whose expansion in the universal
character basis, has polynomial coefficients in $q$ with
non-negative integer coefficients. We conjecture that these
polynomials are one-dimensional sums associated with crystal bases
of finite-dimensional modules over quantized affine algebras for
all nonexceptional affine types. These polynomials satisfy a
Macdonald-type duality.
\end{abstract}

\maketitle

%%% comment this out for distribution
\tableofcontents

\section{Introduction}

\subsection{Universal characters of classical type}

It is well-known that the ring $\La$ of symmetric functions is the
universal character ring of type $A$. That is, for every
$n\in\Z_{>0}$ there is a ring epimorphism $\La\rightarrow
R(GL(n))$ from $\La$ onto the ring of polynomial representations
of $GL(n)$, which sends the Schur function $s_\la$ to the
isomorphism class of the irreducible $GL(n)$-module of highest
weight $\la$.

The ring $\La$ is also the universal character ring for types B,
C, and D. Using identities of Littlewood \cite{L}, Koike and
Terada constructed two distinguished bases $\{s^\vdom_\la\}$ and
$\{s^\hdom_\la\}$ of $\La$ corresponding to the irreducible
characters of the symplectic and orthogonal groups. These two
bases have the same structure constants $d_{\la\mu\nu}$.
%There is
%a third basis $\{s^\cell_\la\}$ of $\La$ that also has structure
%constants $d_{\la\mu\nu}$ which was not considered in \cite{KT}.
%The basis $\{s^\cell_\la \}$ is implicitly defined in \cite{Kl},
%where it was also shown that up to a natural constraint involving
%Schur function expansions, the only bases of $\La$ with structure
%constants $d_{\la\mu\nu}$ are the three bases $\{s^\se_\la\}$ for
%$\se\in\{\cell\,,\hdom\,,\vdom\,\}$.

There is a third basis $\{s^\cell_\la \}$ of $\La$ that also has
structure constants $d_{\la\mu\nu}$ which appears in \cite{Ko}.
This basis is implicitly defined
in \cite{Kl} where it was also shown that up to a natural constraint
involving Schur function expansions, the only bases of $\La$ with
structure constants $d_{\la\mu\nu}$ are the three bases $\{s^\se_\la\}$
for $\se\in\{\cell\,,\hdom\,,\vdom\,\}$.

\subsection{Creation operators}
Bernstein's creation operator $B_r$ \cite{Ze} is a linear
endomorphism of $\La$. The operators $B_r$ create the Schur basis
in the sense that $B_{\la_1} B_{\la_2} \dotsm B_{\la_k} 1 = s_\la$
where $\la=(\la_1,\dots,\la_k)$. %this only works for $\la_k\geq0$
% s_\nu[Z] = J(z^\rho)^{-1} J( z^{\nu+\rho}) for any \nu \in \Z^k
The primary purpose of this
paper is to study the creation operators $B_r^\se$ for the bases
$\{s^\se_\la\}$ and their $q$-analogues $\tB_r^\se$. The
$q$-analogues are defined by the general construction of the
$q$-analogue of any symmetric function operator \cite{Z}.

\subsection{$K$ polynomials}
Jing \cite{J} showed that the type $A$ operators $\tB_r$ create
the modified Hall-Littlewood symmetric functions, whose Schur
function expansion coefficients are the Kostka-Foulkes
polynomials. More generally, there is a parabolic analogue
$\tB_\nu$ of $\tB_r$ which yields the generalized Kostka
polynomials \cite{SZ}.

We consider the $\se$-analogues of these operators. The operators
$\tB_r^\se$ create symmetric functions, which, when expanded in
the basis $s^\se_\la$, yield coefficient polynomials
$d_{\la\mu}(q)$ with nonnegative integer coefficients. These
polynomials are the same for $\se\in\{\cell\,,\hdom\,,\vdom\,\}$
and are $q$-analogues of the multiplicity of the irreducible
character $\chi^\la$ in the tensor product $\bigotimes_i V(\mu_i
\omega_1)$ of type $C_n$ irreducibles indexed by multiples of the
first fundamental weight.

One may consider the straightforward parabolic analogue
$B^\se_\nu$ of the creation operator $B^\se_r$ and its
$q$-analogue $\tB^\se_\nu$. The operators $\tB^\se_\nu$ yield
polynomials that fail to have nonnegative coefficients.

Motivated by considerations detailed below, for each
$\se\in\{\cell\,,\hdom\,,\vdom\,\}$ we define an $\se$-variant
$H^\se_\nu$ of the type $A$ operator $\tB_\nu$. These operators
create symmetric functions whose expansion into the basis
$s^\se_\la$, has coefficients $K^\se\in\N[q]$. The polynomials
$K^\se$ include the polynomials $d_{\la\mu}(q)$ as a special case
and we show that all $K^\se$ are linear combinations of
the $d(q)$ coefficients.

The creation operators give an explicit formula for $K^\se$ in
terms of the generalized Kostka polynomials. The generalized
Kostka polynomials satisfy a Macdonald-type duality \cite{ScWa}
\cite{Sh2}. This implies a duality for the polynomials $K^\se$
which, as a special case, relates the polynomials $K^\vdom$ and
$K^\hdom$.

\subsection{The $X=M=K$ conjecture}
Let $\gggg$ be a Kac-Moody Lie algebra of nonexceptional affine
type, with derived subalgebra $g'$ and simple Lie subalgebra
$\gb$. Let $U_q(\gggg)\supset U'_q(\gggg)\supset U_q(\gb)$ be the
corresponding quantized universal enveloping algebras. The papers
\cite{HKOTT} \cite{HKOTY} assert the existence of a certain family
of finite-dimensional $U'_q(\gggg)$-modules $W^{(k)}_s$ called
Kirillov-Reshetikhin (KR) modules. The KR modules are conjectured
to have crystal bases $B^{k,s}$. It is expected that all
irreducible finite-dimensional $U'_q(\gggg)$-modules that have a
crystal base, are tensor products of KR-modules. Such modules have
the structure of a graded $U_q(\gb)$-module. Their graded
multiplicities are called one-dimensional sums $X$. The above
authors also define the fermionic formula $M$, whose form is
suggested by the Bethe Ansatz, and conjecture that $X=M$.

As the rank of $\gggg$ goes to infinity, the fermionic formula $M$
stabilizes. We call these polynomials the stable fermionic
formulae. Of all the infinite families of affine root systems with
distinguished $0$ node, there are only four different families of
stable fermionic formulae. These four families $M^\se$ for
$\se\in\{\nn\,,\cell\,,\hdom\,,\vdom\,\}$ correspond to the four
bases given by the Schur functions $s^\nn_\la=s_\la$ and the bases
$s^\se_\la$. The family $\se$ associated with a given affine root
system, is determined by the part of the Dynkin diagram near the
$0$ node.

We conjecture that $X=M=K$ for the stable formulae (large rank
case). In type $A^{(1)}_n$ this follows by combining the papers
\cite{KSS} \cite{Sh} \cite{SZ}. It is also known to hold for
$q=1$.
%%% new sentence here
In a separate publication \cite{Sh4} it will be shown that $X=K$
holds for types $\se\in\{\cell\,,\hdom\,\}$ for tensor factors of
the form $B^{1,s}$ using the virtual crystal theory of \cite{OSS}
and a generalization of the Schensted bijection for oscillating
tableaux due to Delest, Dulucq, and Favreau \cite{DDF}.

\subsection{Dual bases and affine Kazhdan-Lusztig polynomials}
The dual bases to $\{ s_\la^\se \}$ for $\se \in \{ \nn\,,
\cell\,, \vdom\,, \hdom \}$ with respect to the standard
scalar product are families
which are elements of the completion of the ring of symmetric
functions.  They are defined as the the bases
$\{ s_\la^{\se\ast} \}$ with the property
$\langle s_\mu^\se, s_\la^{\se\ast}\rangle = \delta_{\la\mu}$.
This implies that $s_\la^{\se\ast}$ for $\se = \nn\,, \cell\,, \vdom\,,
\hdom$ (respectively) is a Schur function
times $\prod_{\alpha \in {\bar \Phi}^+} (1-e^\alpha)^{-1}$ 
(with $x_i = e^{\epsilon_i}$)
where ${\bar \Phi}^+$ are the set of positive roots of type
$A$, $B$, $C$, $D$ (respectively) which
are not positive roots of type $A$.

In an interesting connection we find that 
creation operators for the $s_\la^{\se\ast}$
bases can also be used to define linearly independent sets of elements
related to the Weyl groups of type $A_n$, $B_n$, $C_n$ and $D_n$
which have Lusztig's weight multiplicity polynomials as change
of basis coefficients.
These results will be explored further in a separate paper.

\subsection{Outline of Paper}
We begin in Section \ref{sec:symfunc} with the derivation of
symmetric function identities using the `plethystic' notation \`a
la Garsia. Section \ref{sec:rschars} introduces notation for
classical root systems and translates expressions involving root
systems into the language of symmetric functions.

In section \ref{sec:bosf} we give a plethystic definition for the
bases $\{s^\se_\la\}$ based on Littlewood's formulae.
The bases $\{ s^\se_\la \}$
for $\se \in \{ \nn\,, \vdom\,, \hdom\,, \cell \}$
are defined as special cases of
a single formula involving skewing (dual to multiplication)
operators acting on Schur functions allowing a unified
account of formulas in the rest of the paper.  An important
advantage of our definition, is that any homomorphism on symmetric
functions which transform a basis of one type to another is
easily realized as a substitution on the alphabet of the symmetric
function.

Such algebra homomorphisms allow us to give explicit
formulae for creation operators $B_r^\se$ for the bases
$\{s^\se_\la\}$ in section \ref{sec:bo}. These creation operators
are used to derive several Jacobi-Trudi like identities for the
distinguished bases, in particular recovering some identities of
Weyl and the determinantal definition of the bases
$\{s^\vdom_\la\}$ and $\{s^\hdom_\la\}$ given in \cite{KT}. We
also derive creation operators for the dual bases to $s^\se_\la$.

In section \ref{sec:qanalogues} we study the straightforward
$q$-analogues $\tB^\se_\nu$ of the parabolic creation operators
$B^\se_\nu$. Section \ref{sec:Hse} gives the definition of the $K$
polynomials via the $\se$-analogues $H^\se_\nu$ of the parabolic
Hall-Littlewood creation operator $\tB_\nu$.

Section \ref{sec:XMK} states the $X=M=K$ conjecture.

\section{Plethystic formulae}
\label{sec:symfunc} Let $\La$ be the ring of symmetric functions,
to which we apply the `plethystic notation.' We refer the reader to
\cite{Mac} for symmetric function identities and we will list all
the necessary additional identities we will use in this section but
we will not include a complete exposition on this notation.

Assume that the letters
$X, Y, Z$ and $W$ represent sums of monomials with coefficient $1$ so
that $XY$ represents a product of these sums of monomials.
Expressions like $x \in X$ refers to $x$ being a single
monomial in the multiset of this expression.

$\La$ has a scalar product $\inner{\cdot}{\cdot}$ with respect to
which the Schur functions $\{s_\la\}$ are an orthonormal basis.
The reproducing kernel for this scalar product is $\Omega[XY]$,
where
\begin{equation}
\label{eq:omega} \Omega[X-Y] = \dfrac{\prod_{y\in Y}
(1-y)}{\prod_{x\in X} (1-x)} = \left( \sum_{r \geq 0} (-1)^r s_{(1^r)}[Y] \right)
 \left( \sum_{r \geq 0} s_r[X] \right).
\end{equation}
Cauchy's formula expands the kernel in terms of Schur functions:
\begin{equation} \label{eq:dualbases}
\Omega[XY] = \sum_\la s_\la[X] s_\la[Y].
\end{equation}

Given $P[X]\in\La$, the skewing operator $P[X]^\perp\in\End(\La)$
is the linear operator that is adjoint to multiplication by $P[X]$
with respect to the scalar product. In other words, for all
$P[X]\in\La$,
\begin{equation} \label{eq:skewkern}
  P[X]^\perp(\Omega[XY]) = \Omega[XY] P[Y].
\end{equation}
Taking $P[X]=\Omega[XZ]$ and skewing in the $X$ variables, we have
\begin{equation} \label{eq:kernskewkern}
  \Omega[XZ]^\perp(\Omega[XY]) = \Omega[XY] \Omega[YZ] =
  \Omega[(X+Z)Y].
\end{equation}
\begin{remark}
All skewing operations in this paper are done with respect to symmetric
functions in the $X$ variables only.  Expressions which contain other sets
of variables are considered as coefficients of the linear operator.
In addition, consider $V, W \in End(\Lambda)$ that
do not involve the arbitrary set of variables $Y$, if
$V( \Omega[XY]) = W(\Omega[XY])$, then by (\ref{eq:dualbases}) and by taking
the coefficient of $s_\la[Y]$ it holds that $V( s_\lambda[X]) = W(s_\lambda[X])$.
Since $s_\lambda[X]$
is a basis for $\Lambda$, then $V=W$ as elements of $\End(\Lambda)$.
\end{remark}
Due to the previous remark we have that for all $P[X]\in \La$,
\begin{equation} \label{eq:kernskew}
  \Omega[XZ]^\perp(P[X]) = P[X+Z].
\end{equation}
We compute the commutation of the multiplication operator
$\Omega[XZ]$ and the skewing operator $\Omega[WX]^\perp$, letting
the composite operators act on $\Omega[XY]$.
\begin{equation*}
\begin{split}
  \Omega[XW]^\perp (\Omega[XZ] \Omega[XY])
  &= \Omega[XW]^\perp(\Omega[X(Y+Z)]) \\
  &= \Omega[(X+W)(Y+Z)] \\
  &= \Omega[(X+W)Z] \Omega[XW]^\perp(\Omega[XY]).
\end{split}
\end{equation*}
This yields the operator identity
\begin{equation} \label{eq:kernskewmult}
  \Omega[XW]^\perp \circ \Omega[XZ] =
  \Omega[(X+W)Z] \circ \Omega[XW]^\perp.
\end{equation}
In particular, for all $P[X]\in\La$,
\begin{equation} \label{eq:skewmult}
  \Omega[XW]^\perp \circ P[X] = P[X+W] \circ \Omega[XW]^\perp.
\end{equation}

The comultiplication map $\Delta:\La\rightarrow\La\otimes\La$ may be computed
as follows. Let $P\in \La$, expand $P[X+Y]$ as a sum of products of the form
$P_1[X] P_2[Y]$: $P[X+Y] = \sum_{(P)} P_1[X]P_2[Y]$. Then
$\Delta(P)=\sum_{(P)} P_1 \otimes P_2$.

The adjoint operators $P^\perp$ act on products by \cite[Ex.
I.5.25(d)]{Mac}
\begin{equation} \label{eq:skewprod}
 P^\perp(QR) = \sum_{(P)} P_1^\perp(Q) P_2^\perp(R).
\end{equation}
Given any operator $V\in \End(\La)$, one of the authors \cite{Z}
defined its $t$-analogue $\Vt\in\End(\La)$ by
\begin{equation}\label{tanal}
  \Vt(P[X]) = V^Y(P[tX+(1-t)Y])|_{Y\rightarrow X}
\end{equation}
where $V^Y$ acts on the $Y$ variables and $Y\rightarrow X$ is the
substitution map.

We now apply the above construction for
$V=\Omega[XZ]\circ \Omega[XW]^\perp$. We have that
for $P[X] \in \La$,
\begin{equation*}
  \Omega[YZ] \Omega[YW]^{\perp_Y} P[tX+(1-t)Y] = \Omega[YZ]
  P[tX+(1-t)(Y+W)].
\end{equation*}
Then
\begin{equation*}
\begin{split}
  \Vt(P[X]) &= \Omega[XZ] P[tX+(1-t)(X+W)] \\
  &= \Omega[XZ] P[X+(1-t)W] = \Omega[XZ] \Omega[XW(1-t)]^\perp P[X].
\end{split}
\end{equation*}
By linearity, for all $P[X],Q[X]\in\La$, if $V=P[X] \circ
Q[X]^\perp$, then
\begin{equation} \label{eq:analogmultskew}
  \Vt = P[X] \circ Q[X(1-t)]^\perp.
\end{equation}
For this $V$, at $t=0$ the operator $V$ is recovered:
\begin{equation} \label{eq:analoguezero}
  \Vt|_{t=0} = P[X] \circ Q[X]^\perp = V.
\end{equation}
At $t=1$, the operator
\begin{equation} \label{eq:analogueone}
  \Vt|_{t=1} = P[X] Q[0]
\end{equation}
is multiplication by $P[X] Q[0]$.

For later use we compute the commutation of the operator
$\Omega[s_{(1^2)}[X]]^\perp$ with a multiplication operator.

\begin{prop} \label{pp:e2perpcomm}
\begin{equation} \label{eq:e2perp}
 \Omega[W s_{(1^2)}[X]]^\perp \circ \Omega[ZX] =
 \Omega[ZX+W s_{(1^2)}[Z]] \Omega[W(s_{(1^2)}[X]+ZX)]^\perp.
\end{equation}
\end{prop}
\begin{proof}
\begin{equation*}
\begin{split}
  & \,\,\quad\Omega[W s_{(1^2)}[X]]^\perp ( \Omega[ZX] \Omega[XY])
\\
  &= \Omega[W s_{(1^2)}[X]]^\perp (\Omega[(Y+Z)X]) \\
  &= \Omega[(Y+Z)X] \Omega[W s_{(1^2)}[Y+Z]] \\
  &= \Omega[(Y+Z)X] \Omega[W(s_{(1^2)}[Y]+YZ+s_{(1^2)}[Z])] \\
  &= \Omega[ZX+W s_{(1^2)}[Z]] \Omega[XY] \Omega[W(s_{(1^2)}[Y]+YZ)] \\
  &= \Omega[ZX+W s_{(1^2)}[Z]] (\Omega[W(s_{(1^2)}[X]+ZX)]^\perp (\Omega[XY])).
\end{split}
\end{equation*}
\end{proof}

\section{Root systems and characters}\label{sec:rschars}
In this section we fix notation for the root systems of classical
type and recall Littlewood's plethystic formula and Weyl's
determinantal formulae for the Weyl denominator.

Consider a root system $\Phi$ of type $A_{n-1}$, $B_n$, $C_n$, or
$D_n$, coming from the Lie group $GL_n$, $O_{2n+1}$, $Sp_{2n}$, or
$O_{2n}$ respectively. We use $GL_n$ instead of $SL_n$ for
convenience.

Let $I=\{1,2,\dotsc,n-1\}$ for type $A_{n-1}$ and
$I=\{1,2,\dotsc,n\}$ for types $B_n,C_n,D_n$. Let $\epsilon_i$ be
the $i$-th standard basis vector of $\Z^n$. A set of simple roots
$\Delta=\{\alpha_i\mid i\in I\}$ can be given by
\begin{equation*}
\Delta=\begin{cases}
\{\alpha_i=\epsilon_i-\epsilon_{i+1} \mid 1\le i\le n-1 \} & \text{for type $A_{n-1}$} \\
\{\alpha_i=\epsilon_i-\epsilon_{i+1} \mid 1\le i\le n-1 \} \cup \{\alpha_n=\epsilon_n \} & \text{for type $B_n$} \\
\{\alpha_i=\epsilon_i-\epsilon_{i+1} \mid 1\le i\le n-1 \} \cup \{\alpha_n=2\epsilon_n \} & \text{for type $C_n$} \\
\{\alpha_i=\epsilon_i-\epsilon_{i+1} \mid 1\le i\le n-1 \} \cup \{\alpha_n=\
\epsilon_{n-1}+\epsilon_n \} & \text{for type $D_n$}
\end{cases}
\end{equation*}
The simple roots form a basis for the real Euclidean space $E$ that they span.
For type $A_{n-1}$, $E=\{v\in\R^n\mid \inner{v}{\epsilon} =0 \}$ where $\epsilon=(1,1,\dotsc,1)\in\R^n$.
For types $B_n,C_n,D_n$, $E=\R^n$.

The simple coroots $\{\ac_i\mid i\in I\}\subset E$ are defined by
$\ac_i = \frac{2}{\inner{\alpha_i}{\alpha_i}}\alpha_i$ for $i\in I$,
where $\inner{\cdot}{\cdot}$ is the standard scalar product on $\R^n$.

The fundamental weights $\{\omega_i\mid i\in I\}\subset E^*$ are defined by
the condition that they be a dual basis to the coroot basis of $E$
with respect to the standard scalar product on $\R^n$. The weight
lattice is defined by $P=\bigoplus_{i\in I} \Z \omega_i \subset E^*$ and
the dominant weights are given by $P^+=\bigoplus_{i\in I} \Z_{\ge0} \omega_i$.
In type $A_{n-1}$, $E^* \cong \R^n/\R\epsilon$. For types $B_n,C_n,D_n$ we
shall make the identification $E^*=E$.

The dominant weights are given by
\begin{equation*}
\begin{cases}
\{\omega_i=(1^i,0^{n-i}) \mod \R \epsilon \mid 1\le i\le n-1 \} & \text{for type $A_{n-1}$} \\
\{\omega_i=(1^i,0^{n-i}) \mid 1\le i\le n-1 \}\cup
\{\omega_n=(1/2) \epsilon \} & \text{for type $B_n$} \\
\{\omega_i=(1^i,0^{n-i}) \mid 1\le i\le n-1 \}\cup
\{\omega_n= \epsilon \} & \text{for type $C_n$} \\
\{\omega_i=(1^i,0^{n-i}) \mid 1\le i\le n-2 \}\cup \\
\qquad\{\omega_{n-1}=((1/2)^{n-1},-1/2),
\omega_n=(1/2) \epsilon \} & \text{for type $D_n$}
\end{cases}
\end{equation*}

For $i\in I$ let $s_i\in GL(E)$ be the linear map given by the
reflection through the hyperplane orthogonal to $\alpha_i$:
$s_i(\la)=\la-\inner{\la}{\ac_i}\alpha_i$ for $i\in I$ and $\la\in
E$. The Weyl group $W$ is the group generated by $\{s_i \mid i\in
I\}$. For type $A_{n-1}$, $W$ is the symmetric group $S_n$ on $n$
symbols. For types $B_n$ and $C_n$, $W$ is the hyperoctahedral
group of signed permutations on $n$ symbols. For type $D_n$ $W$ is
the subgroup of the hyperoctahedral group on $n$ symbols
consisting of the signed permutations with an even number of
negative signs.

The root system $\Phi$ is the $W$-orbit of $\Delta$ in $E$.
The set $\Phi^+$ of positive roots is defined by
\begin{equation*}
  \Phi^+ = \{\alpha\in \Phi \mid \inner{\alpha}{\omega_i} \ge 0 \text{ for all $i\in I$. } \}
\end{equation*}
Explicitly
\begin{equation*}
  \Phi^+ = \begin{cases}
  \{ \epsilon_i-\epsilon_j \mid 1\le i<j\le n \} & \text{for type $A_{n-1}$} \\
  \{ \epsilon_i\pm\epsilon_j \mid 1\le i<j\le n \}\cup\{\epsilon_i\mid 1\le i\le n\} &
  \text{for type $B_n$} \\
  \{ \epsilon_i\pm\epsilon_j \mid 1\le i<j\le n \}\cup\{2\epsilon_i\mid 1\le i\le n \} &
  \text{for type $C_n$} \\
  \{ \epsilon_i\pm\epsilon_j \mid 1\le i<j\le n \}&
  \text{for type $D_n$}
  \end{cases}
\end{equation*}
Let $\rho=\frac{1}{2}\sum_{\alpha\in \Phi^+} \alpha=\sum_{i\in I} \omega_i $.
Explicitly
\begin{equation*}
\rho=\begin{cases}
(n-1,n-2,\dotsc,1,0) \mod \R \epsilon & \text{for type $A_{n-1}$} \\
(1/2)(2n-1,2n-3,\dotsc,3,1) & \text{for type $B_n$} \\
(n,n-1,\dotsc,2,1) & \text{for type $C_n$} \\
(n-1,n-2,\dotsc,1,0)& \text{for type $D_n$}
\end{cases}
\end{equation*}

Let $\Z[P]$ be the group $\Z$-algebra of the weight lattice, with
$\Z$-basis $e^\la$ for $\la\in P$. The Weyl group $W$ acts on $P$
and therefore on $\Z[P]$. Let $J=\sum_{w\in W} \minus(w) w$ be the
$W$-antisymmetrization operator on $\Z[P]$. Define the linear map
$\pi$ on $\Z[P]$ by $\pi f = J(e^\rho)^{-1} J(e^\rho f)$.

The Weyl character formula says that the character $\chi^\la$ of the irreducible
module of highest weight $\la\in P^+$, is given by
\begin{equation}
  \chi^\la = \pi(e^\la).
\end{equation}
We define this also for any $\la\in P$. For any $w\in W$ we have
\begin{equation}
  \chi^\la = \minus(w) \chi^{w(\la+\rho)-\rho}.
\end{equation}
It follows that $\chi^\la\not=0$ if and only if $\la+\rho$ has trivial stabilizer
in $W$. In this case there is a unique $w\in W$ such that $w(\la+\rho)\in P^+$.
Thus $\chi^\la$ is either zero or an irreducible character, up to sign.

The Weyl denominator (up to the factor $e^{-\rho}$) is defined by
\begin{equation} \label{eq:Weyldenom}
  \prod_{\alpha\in \Phi^+} (1-e^{-\alpha}) = e^{-\rho} J(e^\rho)
\end{equation}
Let $z_i=e^{\epsilon_i}$ for $1\le i\le n$,
$Z=(z_1,z_2,\dotsc,z_n)$, $z_i^*=z_i^{-1}$, and
$Z^*=z_1^*+z_2^*+\dotsm+z_n^*$. Then the Weyl denominator for these four cases
is given in our notation as
\begin{equation}
\begin{split}
R^A(Z) &= R(Z) = \prod_{1\le i<j\le n} (1-z_jz_i^*) \\
R^B(Z) &= R(Z) \Omega[-f_\cell\,[Z^*]] \\
R^C(Z) &= R(Z) \Omega[-f_\hdom\,[Z^*]] \\
R^D(Z) &= R(Z) \Omega[-f_\vdom\,[Z^*]]
\end{split}
\end{equation}
where we have defined
\begin{equation} \label{eq:kern}
\begin{split}
\KK_\nn[X] &= 0 \\
\KK_\cell[X] &= s_1[X]+s_{(1^2)}[X] \\
\KK_\hdom[X] &= s_2[X] \\
\KK_\vdom[X] &= s_{(1^2)}[X]
\end{split}.
\end{equation}
We remark that due to symmetric function identies, we know that
$s_1[Z^\ast] = Z^\ast$, $s_{(1^2)}[Z^\ast] = \sum_{i<j} z_i^\ast z_j^\ast$
and $s_2[Z^\ast] = \sum_{i\leq j} z_i^\ast z_j^\ast$.  Along with equation
(\ref{eq:omega}), we may arrive at other expressions for $R^X(Z)$, but for
manipulations we have chosen to express them as symmetric function identities.

The following is essentially due to Weyl.

\begin{prop} \label{pp:denomdet} Consider the algebra automorphism of $\Z[P]$
induced by $z_i\mapsto z_{n+1-i}^*$. Letting $\Rb$ denote the image of the
Weyl denominator $R$ under this map, we have the $n\times n$ determinantal
formulae
\begin{align*} \label{eq:Weyldet}
  R(Z) = \Rb^A(Z) &= \det | z_i^{i-j} |  \\
  R(Z) \Omega[-\KK_\cell[Z]] = \Rb^B(Z) &= \det | z_i^{i-j}-z_i^{i+j-1} | \\
  R(Z) \Omega[-\KK_\hdom[Z]] = \Rb^C(Z) &= \det | z_i^{i-j}-z_i^{i+j} | \\
  R(Z) \Omega[-\KK_\vdom[Z]] = \Rb^D(Z) &= \frac{1}{2}\, \det | z_i^{i-j} + z_i^{i+j-2} |
\end{align*}
\end{prop}
\begin{proof} Let $\rrho$ be the reverse of $\rho$. The algebra automorphism
sends $e^{w\rho-\rho}$ to $e^{\rrho-w\rrho}$. For type $B_n$ we have that
$\rrho=\frac{1}{2}(1,3,\dotsc,2n-1)$ and
\begin{align*}
R(Z) \Omega[-Z_n - s_{(1^2)}[Z_n]] &= \sum_{w \in S_n} \epsilon(w)
\sum_{b \in \{ \pm 1\}^n} \left( \prod_{i=1}^n b_i \right)
z^{\rrho - \alpha w \rrho}\\
&= \sum_{w \in S_n} \epsilon(w) \sum_{b \in \{ \pm 1\}^n}
\prod_{i=1}^n b_i z_i^{ i- 1/2 - b_i (w_i-1/2)}\\
&= \sum_{w \in S_n} \epsilon(w) \prod_{i=1}^n
\left(z_i^{i-1/2-w_i+1/2} - z_i^{i-1/2+w_i-1/2} \right)\\
&= \det \left| z_i^{i-j} - z_i^{i+j-1} \right|_{1 \leq i,j \leq n}
\end{align*}
The other types are similar.
\end{proof}

\section{Bases of symmetric functions}\label{sec:bosf}
We consider four bases of the algebra $\La$ of symmetric
functions. The first is the basis $\{s_\la\}$ of Schur functions.

For $\se\in\{\,\nn, \cell\,,\,\vdom\,,\,\hdom\,\}$, define the
symmetric function
\begin{equation}
\label{eq:basis} s^\se_\la[X] = \Omega[-\KK_\se]^\perp s_\la[X]
\end{equation}
with $\KK_{\se}$ as in \eqref{eq:kern}. All of the families
$\{s^\se_\la \}$ are bases of $\La$, due to the inverse formula
\begin{equation}
\label{eq:inv} s_\la[X] = \Omega[\KK_\se]^\perp s^\se_\la[X].
\end{equation}
Of course $s^\nn_\la=s_\la$ is the basis of Schur functions, which
are the universal characters for the special/general linear
groups. The bases $\{\ssp_\la\}$ and $\{\so_\la\}$ appear in
\cite{KT} as the universal characters for the symplectic and
orthogonal groups respectively; see \eqref{eq:charse}. The basis
$\{\sd_\la\}$ is not mentioned in \cite{KT} but appears implicitly
in \cite{Kl}.

To explain the notation for the bases, let $\PPP^\se$ be the set
of partitions that can be tiled using the shape $\se$ without
changing the orientation of the tile $\se$. In other words,
$\PPP^\nn=\{\nn\}$ is the singleton set containing the empty
partition, $\PPP^\cell$ is the set of all partitions, $\PPP^\hdom$
is the set of partitions with even rows, and $\PPP^\vdom$ is the
set of partitions with even columns.

Littlewood's formulae give the Schur function expansions of
$\Omega[\pm \KK_\se]$.
\begin{align}
\label{eq:skewpos} \Omega[\KK_\se] &= \sum_{\la\in P^\se}
s_\la[X] \\
\label{eq:ivdom} \Omega[-\KK_\vdom\,] &=
\sum_{\mu=(\alpha_1,\dotsc,
\alpha_p|\alpha_1-1,\dotsc,\alpha_p-1)} (-1)^{|\mu|/2} s_\mu[X]
= \prod_{i<j} 1- x_i x_j
\\
\label{eq:ihdom} \Omega[-\KK_\hdom\,] &= \sum_{\mu=(\alpha_1-1,
\ldots, \alpha_p-1|\alpha_1,\ldots,\alpha_p)} (-1)^{|\mu|/2}
s_\mu[X] = \prod_{i\leq j} 1- x_i x_j
\\
\label{eq:icell} \Omega[-\KK_\cell\,] &= \sum_{\mu=\mu^t}
(-1)^{(|\mu|+d(\mu))/2} s_\mu[X] = \prod_{i<j} 1- x_i x_j \prod_i 1- x_i,
\end{align}
where $(\alpha|\beta)$ is Frobenius' notation for a partition and
$d(\mu)$ is the size of largest diagonal.

We introduce notation for the linear maps that change among the
bases $\{s^\se_\la\}$ for various $\se$. In plethystic formulae
let $\minus$ represent a variable that has been specialized to the
scalar $-1$.  We will consider $\minus$ a special element
with the property $\minus^2=1$ and
\begin{equation}
  \Omega[\minus X - \minus Y] = \dfrac{\prod_{y\in Y}
(1+y)}{\prod_{x\in X} (1+x)}
\end{equation}
 For $\se,\heartsuit\in
\{\nn,\cell\,,\vdom\,,\hdom\,\}$ define the linear isomorphism
$i_\se^\heartsuit:\La\rightarrow\La$ by
\begin{equation} \label{eq:changebasis}
  i_\se^\heartsuit(s^\se_\la[X]) = s^\heartsuit_\la[X]
\end{equation}
for all $\la$. It is given by
\begin{equation}
  i_\se^\heartsuit = \Omega[\KK_\se-\KK_\heartsuit]^\perp.
\end{equation}

\begin{prop} \label{pp:dchange} For all $P\in \La$,
\begin{alignat}{3}
i_\vdom^\cell P[X] &= P[X-1] \qquad & i_\cell^\vdom P[X] &= P[X+1]
\\
i_\vdom^\hdom P[X] &= P[X-1-\minus] \qquad & i_\hdom^\vdom
P[X] &= P[X+1+\minus] \\
i_\cell^\hdom P[X] &= P[X-\minus] \qquad & i_\hdom^\cell P[X] &=
P[X+\minus]
\end{alignat}
\end{prop}
\begin{proof} We have
\begin{align}
\label{eq:vdomcell}
  \Omega[\KK_\vdom-\KK_\cell] &= \Omega[-X] \\
\label{eq:vdomhdom}
  \Omega[\KK_\vdom-\KK_\hdom] &= \Omega[-p_2[X]] =
\prod_{x\in X} (1-x^2) =
  \Omega[-(1+\minus)X]
\end{align}
All of the formulae follow from these and \eqref{eq:kernskew}.
\end{proof}

In particular, since substitution maps are algebra homomorphisms,
one has the following result, which was obtained in \cite{KT} for
the pair $\vdom\,$ and $\hdom\,$.

\begin{cor} \label{cor:algiso} $i_\se^\heartsuit$ is an algebra
isomorphism for $\se,\heartsuit\in\{\cell\,,\vdom\,,\hdom\,\}$.
\end{cor}

Define the structure constants ${}^\se c^\la_{\mu\nu}$ by
\begin{equation} \label{eq:spc}
s^\se_\mu[X] s^\se_\nu[X] = \sum_\la {}^\se c^\la_{\mu\nu}
s^\se_\la[X].
\end{equation}
The coefficient ${}^\nn c^\la_{\mu\nu}$ is the ordinary
Littlewood-Richardson coefficient $c^\la_{\mu\nu}$. By Corollary
\ref{cor:algiso}, the other three sets of structure constants
coincide (this is proved in \cite{KT} for $\vdom$ and $\hdom$);
call this common structure constant $d_{\la\mu\nu}$.

Kleber \cite{Kl} showed that the three bases $\{s^\se_\la\}$ for
$\se\in\{\cell\,,\vdom\,,\hdom\,\}$ are the unique bases with this
property, assuming a certain kind of expansion in terms of Schur functions.

\begin{theorem} \label{thm:only} \cite{Kl}
Suppose $\{v_\la\}$ is a basis of the ring of symmetric functions
such that
\begin{equation}
  v_\mu v_\nu = \sum_\la d_{\la\mu\nu} v_\la
\end{equation}
for all $\mu,\nu$ and that
\begin{equation}
  s_\la \in v_\la + \sum_{\mu<\la} \NN \,v_\mu
\end{equation}
where $\mu\le \la$ means that $\mu_1+\dotsm+\mu_i\le
\la_1+\dotsm+\la_i$ for all $i$ (but $\mu$ and $\la$ need not have
the same number of cells). Then $\{v_\la\}$ must be one of the
bases $\{\ssp_\la\}$, $\{\so_\la\}$, or $\{\sd_\la\}$.
\end{theorem}

\begin{example} The element $s^\se_{(433)}$ is expanded in the Schur
basis for $\se\in \{\cell\,,\hdom\,,\vdom\,\}$. Each Schur
function is represented by the diagram for the partition indexing it.
These expansions may computed by \eqref{eq:basis}, \eqref{eq:ivdom},
\eqref{eq:ihdom}, \eqref{eq:icell} and the Littlewood-Richardson
rule.

\begin{equation}
\begin{split}
s_{(433)}^\hdom = &~\young{&&\cr&&\cr&&&\cr} -
\young{\cr&&\cr&&&\cr}-\young{&\cr&&\cr&&\cr}+\young{&\cr&&&\cr}+
\young{&&\cr&&\cr}+\young{\cr&\cr&&\cr}\\
&-\young{&&&\cr}-\young{\cr&&\cr}
-\young{&\cr&\cr}+\young{&\cr}\\
s_{(433)}^\vdom = &~\young{&&\cr&&\cr&&&\cr}
-\young{&\cr&\cr&&&\cr}-\young{&\cr&&\cr&&\cr}
+\young{\cr&\cr&&\cr}+\young{&\cr&\cr&\cr}-\young{\cr\cr&\cr}\\
s_{(433)}^\cell =&~\young{&&\cr&&\cr&&&\cr}
-\young{&\cr&&\cr&&&\cr}-\young{&&\cr&&\cr&&\cr}
+\young{\cr&\cr&&&\cr}+\young{\cr&&\cr&&\cr}+\young{&\cr&\cr&&\cr}
-\young{\cr\cr&&&\cr}-\young{\cr&\cr&&\cr}\\
&-\young{&\cr&&\cr}-\young{\cr&\cr&\cr}
+\young{\cr&&\cr}+\young{&\cr&\cr}+\young{\cr\cr&\cr}
-\young{&\cr}-\young{\cr\cr}+\young{\cr}
\end{split}
\end{equation}
\end{example}

The structure constants $d_{\la\mu\nu}$ can be expressed in terms
of the Littlewood-Richardson coefficients $c^\la_{\mu\nu}$, using
the Newell-Littlewood formula, which we rederive here.

\begin{prop} \label{pp:dc} \cite{L2} \cite{N}
\begin{equation} \label{eq:dc}
  d_{\la\mu\nu} = \sum_{\rho,\sigma,\tau\in\Par} c^\mu_{\rho\tau} c^\nu_{\sigma\tau} c^\la_{\rho\sigma}.
\end{equation}
\end{prop}
\begin{proof} Note that $\Omega[s_{(1^2)}[X+Y]]=\Omega[s_{(1^2)}[X]]\Omega[XY]
\Omega[s_{(1^2)}[Y]]$.
By \eqref{eq:skewprod} we have
\begin{equation}
\begin{split}
  d_{\la\mu\nu} &= \inner{s^\vdom_\mu s^\vdom_\nu}{s^\vdom_\la}
  = \inner{\Omega[s_{(1^2)}[X]]^\perp\left(s^\vdom_\mu s^\vdom_\nu\right)}{s_\la} \\
  &= \sum_\tau \inner{(\Omega[s_{(1^2)}[X]]^\perp s_\tau[X]^\perp s^\vdom_\mu)
   (\Omega[s_{(1^2)}[X]]^\perp s_\tau[X]^\perp s^\vdom_\mu)}{s_\la} \\
   &= \sum_\tau \inner{(s_\tau^\perp s_\mu)(s_\tau^\perp s_\nu)}{s_\la}
   = \sum_\tau \inner{s_{\mu/\tau} s_{\nu/\tau}}{s_\la} \\
   &= \sum_{\rho,\sigma,\tau} c^\mu_{\rho\tau} c^\nu_{\sigma\tau} c^\la_{\rho\sigma}.
\end{split}
\end{equation}
\end{proof}

\begin{cor} $d_{\la\mu\nu}$ is symmetric in its three arguments.
\end{cor}

The next result is an immediate corollary of Proposition \ref{pp:dc} and the
transpose symmetry of Littlewood-Richardson coefficients
\begin{equation} \label{eq:LRtr}
  c^{\la^t}_{\mu^t\nu^t} = c^\la_{\mu\nu}.
\end{equation}
Equation \eqref{eq:LRtr} says that the involution $\omega:\La\rightarrow\La$
given by $\omega(s_\la)=s_{\la^t}$, is an algebra isomorphism.  This implies
the following corollary (see \cite[Theorem 2.3.4]{KT}).

%??? we need the reference to [KT] for this
\begin{cor} \label{cor:dtr} $d_{\la^t\mu^t\nu^t} = d_{\la\mu\nu}$.
\end{cor}

The Schur functions are self dual with respect to the standard
scalar product of the symmetric functions, that is we have
$\left< s_\la[X], s_\mu[X] \right> = \delta_{\la\mu}$.  Since
by definition,
$$\left< \Omega[-f_\se[X]]^\perp f[X], g[X] \right> =
\left< f[X], \Omega[-f_\se[X]] g[X] \right>,$$
the bases $s^\se_\la[X] = \Omega[- f_\se[X]]^\perp s_\la[X]$ will
be dual to the functions $s^{\se\ast}_\la[X] := \Omega[f_\se[X]] s_\la[X]$
since
\begin{align*}
\left< s^{\se\ast}_\la[X], s^\se_\mu[X] \right>
&= \left< \Omega[f_\se[X]] s_\la[X], \Omega[-f_\se[X]]^\perp s_\mu[X] \right>\\
&= \left< s_\la[X], \Omega[f_\se[X]]^\perp \Omega[-f_\se[X]]^\perp s_\mu[X] \right>\\
&= \left< s_\la[X], s_\mu[X] \right> = \delta_{\la\mu}.
\end{align*}

The elements $s^{\se\ast}_\la[X]$ are elements of the
completion of the symmetric functions, $\hat \La$, and
computation of the expansion in the Schur basis can be
computed using equation (\ref{eq:skewpos}).

We shall only consider characters $\chi^\la$ where $\la\in\Z^n\cap
P^+$ using the above realizations of $P$. The connection between
the characters $\chi^\la$ and the symmetric functions $s^\se_\la$
is due to Littlewood. The character
$\chi^\la:G\rightarrow\mathbb{C}$, defined by the trace of the
action of $g\in G$ on the highest weight module $V^\la$ of $G$, is
a polynomial function in the eigenvalues of $g$. A typical element
of the group $G$ has eigenvalues of the form
\begin{equation*}
  \begin{cases}
  (z_1,\dotsc,z_n) & \text{for $G=GL_n$} \\
  (z_1,\dotsc,z_n,1,z_1^{-1},\dotsc,z_n^{-1}) & \text{for $G=O_{2n+1}$} \\
  (z_1,\dotsc,z_n,z_1^{-1},\dotsc,z_n^{-1}) & \text{for $G=Sp_{2n}$ or $O_{2n}$.}
  \end{cases}
\end{equation*}
With $g\in G$ having the above eigenvalues and
$Z= z_1+\cdots+z_n$, we have
\begin{equation} \label{eq:charse}
  \chi^\la(g) = \begin{cases}
  s^\nn_\la[Z] & \text{for $G=GL_n$} \\
  s^\hdom_\la\,[Z+1+Z^*] & \text{for $G=O_{2n+1}$} \\
  s^\vdom_\la\,[Z+Z^*] & \text{for $G=Sp_{2n}$} \\
  s^\hdom_\la\,[Z+Z^*] & \text{for $G=O_{2n}$.}
  \end{cases}
\end{equation}
Because of relations between the bases we also have
$s^\cell_\la[Z+Z^\ast+1] = \chi^\la(g)$ for $G=Sp_{2n}$.
Proctor's odd symplectic group $Sp_{2n+1}$ \cite{P}, which is
neither simple nor reductive, has an indecomposable representation
with trace $\chi^\la$, which restricted to elements $g\in
Sp_{2n+1}$ of determinant $1$, satisfies
$\chi^\la(g)=s^\vdom_\la[Z+1+Z^*]$.

\section{Creation operators for the bases $s^\se_\la$ and determinants}
\label{sec:bo}

\subsection{The Schur basis}
The Schur functions $\{s_\la\mid\la\in\Par\}$ are the unique
family of symmetric functions indexed by partitions, which for
$\la=(r)$ are given by
\begin{equation} \label{eq:onerowSchur}
  s_r[X] = \begin{cases}
    s_r[X] & \text{if $r\in\Z_{\ge0}$} \\
    0 & \text{otherwise,}
  \end{cases}
\end{equation}
and for $\la\in \Par$ are given by the Jacobi-Trudi determinant
\begin{equation}\label{eq:JT}
  s_\la[X] = \det \left| s_{\la_i-i+j}[X]\right|_{1\le i,j\le
  \ell(\la)}
\end{equation}
where $\ell(\la)$ is the number of parts of $\la$. We may also
define $s_\nu[X]$ for $\nu\in\Z^n$ using \eqref{eq:onerowSchur}
and \eqref{eq:JT}. It follows that for all $w$ in the symmetric
group $W=S_n$,
\begin{equation} \label{eq:Schurskewsym}
  s_\nu[X] = \minus(w) s_{w(\nu+\rho)-\rho}[X]
\end{equation}
where $\rho=(n-1,n-2,\dotsc,1,0)\in\Z^n$.

We now recall the creation operators for the Schur basis.
Bernstein's operators $\{B_r\mid r\in\Z\} \subset\End(\La)$ are
defined by the generating function
\begin{equation} \label{eq:berndef}
 B(z) = \sum_{r\in\Z} B_r z^r = \Omega[zX] \Omega[-z^* X]^\perp
\end{equation}
where $z^*=1/z$. For $\nu\in\Z^n$, define
\begin{equation*}
  B_\nu = B_{\nu_1} \circ B_{\nu_2}\circ \dotsm \circ B_{\nu_n} \in\End(\La).
\end{equation*}
We claim that
\begin{equation} \label{eq:crSchur}
  B_\nu 1 = s_\nu[X].
\end{equation}
The proof is included since the proofs for the other bases
$s^\se_\la$ are very similar.

This result follows from a formula for the composition of Bernstein operators.
Let $Z=z_1+z_2+\cdots+z_n$ and
$Z^*=z_1^*+\cdots+z_n^*$. Define
\begin{equation}
B(Z)=\sum_{\nu\in\Z^n} z^\nu B_\nu = B(z_1)\dotsm B(z_n).
\end{equation}
We have
\begin{equation} \label{eq:twobern}
\begin{split}
B(z) B(w) &= \Omega[zX] \Omega[-z^* X]^\perp \Omega[wX]
\Omega[-w^* X]^\perp \\
  &= \Omega[zX] \Omega[w(X-z^*)] \Omega[-z^* X]^\perp \Omega[-w^* X]^\perp \\
  &= \Omega[-wz^*] \Omega[(z+w)X] \Omega[-(z^*+w^*)X]^\perp.
\end{split}
\end{equation}
Iterating \eqref{eq:twobern} yields
\begin{equation} \label{eq:itBern}
  B(Z) = R(Z) \Omega[ZX] \circ \Omega[-Z^* X]^\perp.
\end{equation}
Rewriting $R(Z)$ by Proposition \ref{pp:denomdet} for type $A_{n-1}$,
applying $B(Z)$ to $1\in \La$ and taking the coefficient at $z^\nu$,
we have that
\begin{equation}
\begin{split}
B_\nu 1 &= \Omega[ZX] \det | z_i^{i-j} |_{1\le i,j\le n} \coeff_{z^{\nu}}\\
&= \det| s_{\nu_i - i + j}[X] |_{1 \le i,j \le n} = s_\nu[X].
\end{split}
\end{equation}
\begin{comment}
\begin{equation*}
  \sum_{\nu\in P^+} \chi^\nu(Z) B_\nu = \Omega[ZX]\Omega[-Z^*X]^\perp.
\end{equation*}
\end{comment}

\subsection{Creating the bases $s^\se$}
Define the operators $B^\se_\nu\in\End(\La)$ for $\nu\in\Z^n$ by
\begin{equation} \label{eq:dberndef}
  B^\se(Z) = \sum_{\nu\in\Z^n} z^\nu B^\se_\nu =
 i_\nn^\se \circ B(Z) \circ i_\se^\nn.
\end{equation}
For $\nu\in\Z^n$ it follows from \eqref{eq:crSchur} and \eqref{eq:changebasis} that
\begin{equation} \label{eq:dcreaterows}
  B_{\nu_1}^\se \dotsm B_{\nu_n}^\se 1 = B^\se_{\nu} 1 =
  s^\se_\nu[X].
\end{equation}
The operator $B^\se(Z)$ has the following plethystic formula.

\begin{prop} \label{pp:dbern} For
$\se\in\{\cell\,,\vdom\,,\hdom\,\}$,
\begin{equation} \label{eq:dbern}
  B^\se(Z) = R(Z) \Omega[-\KK_\se[Z]] \Omega[ZX]
\Omega[-(Z+Z^*)X]^\perp.
\end{equation}
\end{prop}
\begin{proof} First let $\se=\vdom\,$.
By \eqref{eq:itBern} and \eqref{eq:e2perp} with $W=-1$ (as an
alphabet), we have
\begin{equation*}
\begin{split}
B^\vdom(Z) &= i_\nn^\vdom B(Z) i_\vdom^\nn \\
&= \Omega[-s_{(1^2)}[X]]^\perp R(Z)
 \Omega[ZX]\Omega[-Z^*X]^\perp \Omega[s_{(1^2)}[X]]^\perp \\
 &= R(Z)\Omega[-s_{(1^2)}[Z]] \Omega[ZX] \Omega[-ZX]^\perp
\Omega[-Z^*X]^\perp \\ &= R(Z) \Omega[-\KK_\vdom\,[Z]] \Omega[ZX]
\Omega[-(Z+Z^*)X]^\perp.
\end{split}
\end{equation*}
Next let $\se=\cell\,$.
\begin{equation*}
\begin{split}
  B^\cell(Z) &= i_\vdom^\cell B^\vdom(Z) i_\cell^\vdom \\
  &= \Omega[-X]^\perp R(Z) \Omega[-\KK_\vdom[Z]] \Omega[ZX]
  \Omega[-(Z+Z^*)X]^\perp \Omega[X]^\perp \\
  &= R(Z) \Omega[-s_{(1^2)}[Z]] \Omega[Z(X-1)]
    \Omega[-(Z+Z^*)X]^\perp \\
  &= R(Z) \Omega[-\KK_\cell[Z]] \Omega[ZX]
\Omega[-(Z+Z^*)X]^\perp.
\end{split}
\end{equation*}
The proof for $\se=\hdom\,$ is similar.
\end{proof}

It follows that
\begin{align}
\label{eq:Bcell}
  B^\cell(Z) &= \Omega[-Z] B^\vdom(Z) \\
\label{eq:Bhdom}
  B^\hdom(Z) &= \Omega[-(1+\minus)Z] B^\vdom(Z)
\end{align}

\subsection{Determinantal formulae}
Recall that the Schur functions satisfy the Jacobi-Trudi identity
\eqref{eq:JT}. The other three bases satisfy a common
determinantal formula due to Weyl for $\ssp$ and
$\so$. See \cite[Thm. 2.3.3]{KT}.

\begin{prop} \label{pp:det} For $\se\in\{ \cell\,,\,\vdom\,,\,\hdom\,\}$ the
basis $\{s^\se_\la \mid \la\in\Par\}$ of $\La$ is characterized by
\begin{equation} \label{eq:onerow}
\begin{split}
s^\vdom_r &= s_r \\
s^\cell_r &= s_r - s_{r-1} \\
s^\hdom_r &= s_r - s_{r-2}
\end{split}
\end{equation}
for $r\in\Z$ and
\begin{equation}
\label{eq:dJT} s^\se_\la = \dfrac{1}{2}
\det \left | s^\se_{\la_i-i+j}+s^\se_{\la_i-i-j+2} \right |_{1\le i,j\le \ell(\la)}
\end{equation}
\end{prop}
\begin{proof} Take $n=\ell(\la)$.
Using the type $D_n$ denominator formula in Proposition \ref{pp:denomdet},
Proposition \ref{pp:dbern} for $\se=\vdom\,$, applying $B^\vdom(Z)$ to $1\in\La$ and
taking the coefficient of $\nu\in\Z^n$ we have that
\begin{equation*}
\begin{split}
  s^\vdom_\nu &= B^\vdom(Z)1\coeff_{z^\nu} \\
  &= \frac{1}{2} \Omega[ZX] \det \left| z_i^{i-j} +  z_i^{i+j-2} \right|\coeff_{z^\nu} \\
  &= \frac{1}{2} \det \left| s_{\nu_i-i+j} + s_{\nu_i+i+j-2}
  \right|.
\end{split}
\end{equation*}
This proves the formula for $\se=\vdom\,$. For $\se\in\{\cell\,,\,\hdom\,\}$
apply the algebra isomorphism $i^\vdom_\se$.
\end{proof}

\begin{example}
\begin{equation*}
s_{(4,3,3)}^\hdom = \frac{1}{2} \left| \begin{array}{ccc}
2 s_4 & s_5 + s_{3} & s_{6} + s_2\\
2 s_2 & s_3 + s_1 & s_4 + 1\\
2 s_1 & s_2 + 1 & s_3
\end{array}
\right|
\end{equation*}
\begin{equation*}
s_{(4,3,3)}^\vdom = \frac{1}{2} \left| \begin{array}{ccc}
2 (s_4-s_2) & s_5 - s_{1} & s_{6}-s_4 + s_2-1\\
2 (s_2-1) & s_3 & s_4-s_2 + 1\\
2 s_1 & s_2 & s_3 - s_1
\end{array}
\right|
\end{equation*}
\begin{equation*}
s_{(4,3,3)}^\cell = \frac{1}{2} \left| \begin{array}{ccc}
2 (s_4-s_3) & s_5-s_4 + s_{3}-s_2 & s_{6}-s_5 + s_2-s_1\\
2 (s_2-s_1) & s_3-s_2 + s_1-1 & s_4-s_3 + 1\\
2 (s_1-1) & s_2-s_1 + 1 & s_3-s_2
\end{array}
\right|
\end{equation*}
\end{example}

Amusingly, if the other denominator identities in Proposition \ref{pp:denomdet} are
used, we obtain more determinantal formulae for $s^\se_\la$. For the denominator
of type $C_n$ we have the following result. For the case $\se=\hdom\,$ see
\cite[Thm. 2.3.3]{KT}.

\begin{prop}
\begin{equation*}
s_\la^{\se} = \det \left| {\tilde s}_{\la-i+j}^\se
- {\tilde s}_{\la_i-i-j}^\se \right|_{1 \leq i,j \leq \ell(\la)}
\end{equation*}
for $\se \in \{ \cell\,, \vdom\,, \hdom\, \}$
where ${\tilde s}^\hdom_r = s_r$, ${\tilde s}^\vdom_r = s_r + s_{r-2} + s_{r-4} + \cdots$, and ${\tilde s}^\cell_r = s_{r} - s_{r-1} + s_{r-2} - \cdots + (-1)^r$.
\end{prop}

\begin{example}
\begin{equation*}
s_{(4,3,3)}^{\hdom} = \det \left| \begin{array}{ccc} s_4 - s_2&s_{5}-s_1&s_{6}-1\\
s_{2}-1 & s_{3} & s_{4}\\
s_{1}&s_{2}&s_{3}
\end{array} \right|
\end{equation*}
\begin{equation*}
s_{(4,3,3)}^{\vdom} = \det \left| \begin{array}{ccc} s_4&s_{5}+s_3&s_{6}+s_4+s_2\\
s_{2} & s_{3}+s_1 & s_{4}+s_2+1\\
s_{1}&s_{2}+1&s_{3}+s_1
\end{array} \right|
\end{equation*}
\begin{equation*}
s_{(4,3,3)}^{\cell} = \det \left| \begin{array}{ccc} s_4-s_3&s_{5}-s_4+s_3-s_2&s_{6}-s_5+s_4-s_3+s_2-s_1\\
s_{2}-s_1 & s_{3}-s_2+s_1-1 & s_{4}-s_3+s_2-s_1+1\\
s_{1}-1&s_{2}-s_1+1&s_{3}-s_2+s_1-1
\end{array} \right|
\end{equation*}
\end{example}

For the type $B_n$ denominator we have new identities.

\begin{prop}
\begin{equation*}
s_\la^{\se} = \det \left| {\overline s}_{\la_i-i+j}^\se
- {\overline s}_{\la_i-i-j+1}^\se \right|_{1 \leq i,j \leq \ell(\la)}
\end{equation*}
for $\se \in \{ \cell\,, \vdom\,, \hdom\, \}$ where ${\overline
s}^\cell_r = s_r$, ${\overline s}^\hdom_r = s_r + s_{r-1}$, and
${\overline s}^\vdom_r = s_{r} + s_{r-1} + s_{r-2} + \cdots + 1$.
\end{prop}

\begin{example}
\begin{equation*}
s_{(4,3,3)}^\cell = \det \left| \begin{array}{ccc} s_4-s_{3}&s_5-s_{2}&s_6-s_1\\
s_2-s_1&s_3-1&s_4\\
s_1-1&s_2&s_3
\end{array} \right|
\end{equation*}
\begin{equation*}
s_{(4,3,3)}^\vdom = \det \left| \begin{array}{ccc} s_4&s_5+s_{4}+s_3&s_6+s_5+s_4+s_3+s_2\\
s_2&s_3+s_2+s_1&s_4+s_3+s_2+s_1+1\\
s_1&s_2+s_1+1&s_3+s_2+s_1+1
\end{array} \right|
\end{equation*}
\begin{equation*}
s_{(4,3,3)}^\hdom = \det \left| \begin{array}{ccc} s_4-s_2&s_5+s_4-s_{2}-s_1&s_6+s_5-s_1-1\\
s_2-1&s_3+s_2-1&s_4+s_3\\
s_1&s_2+s_1&s_3+s_2
\end{array} \right|
\end{equation*}
\end{example}

Finally, we note that all these formulae have analogues involving
the symmetric functions $s_{(1^k)}$ instead of the functions $s_k$.
Recall the algebra isomorphism defined by the
involutive map $\omega:\La\rightarrow\La$ given by
$\omega(s_\la)=s_{\la^t}$ where $\la^t$ is the transpose of the
partition $\la$. Define the map $\omega_\se:\La\rightarrow\La$ by
$i_\nn^\se \circ \omega \circ i_\se^\nn$. By definition
\begin{equation} \label{eq:setranspose}
  \omega_\se(s_\la^\se)=s_{\la^t}^\se.
\end{equation}

We now rewrite \eqref{eq:charse}. Using the definition of
$s^\se_\la$ in the variables $W$, multiplying by $s_\la[X]$,
summing over $\la$, and \eqref{eq:skewkern} we have
\begin{equation*}
\Omega[-f_\se[X]] \Omega[WX] = \sum_{\la\in\Par} s_\la^\se[W] s_\la[X].
\end{equation*}
Multiplying by $\Omega[f_\se[X]]$, we have
\begin{equation} \label{eq:Cauchy}
  \Omega[WX] = \Omega[f_\se[X]] \sum_{\la\in\Par} s_\la^\se[W] s_\la[X].
\end{equation}
Applying this formula with $W=Z,Z+Z^*,Z+1+Z^*$ and
$\se\in\{\nn,\vdom\,,\hdom\,\}$ as in \eqref{eq:charse}, setting
$X$ to be a set of $n$ variables and restricting the sum to
$\la\in\Par^n$, one obtains Littlewood's Cauchy-type formulae
\cite[Lemma 1.5.1]{KT}.

\begin{prop} \label{pp:setranspose}
The map $\omega_\se$ is an algebra isomorphism satisfying
\begin{equation} \label{eq:omegase}
  \omega_\se(s_\la^\se) = \omega(s_\la^{\se^t}).
\end{equation}
where $\nn^t=\nn$, $\vdom\,^t=\hdom\,$, $\hdom\,^t=\vdom\,$, and $\cell\,^t=\cell\,$.
\end{prop}
\begin{proof} $\omega_\se$ is an algebra isomorphism by
\eqref{eq:setranspose} and Corollary \ref{cor:dtr}. To prove
\eqref{eq:omegase} one may reduce to the case $\la=(r)$, since the
$s_r^\se$ are algebra generators of $\La$ by Proposition
\ref{pp:det}. For $\la=(r)$, equation \eqref{eq:omegase} is easy
to check directly using the definitions and the explicit formulae
\eqref{eq:ivdom}, \eqref{eq:ihdom}, and \eqref{eq:icell} for
$\Omega[-f_\se[X]]$.
\end{proof}

By applying the algebra isomorphism \eqref{eq:omegase} to the
formulae in Proposition \ref{pp:det} and using \eqref{eq:omegase}
for $\la=(r)$, one obtains determinantal formulae for
$s_{\la^t}^\se$ in terms of $s^{\se^t}_{1^r}=\omega(s^\se_r)$.
These are simple expressions in the $s_{(1^r)}$ symmetric
functions, due to \eqref{eq:onerow}:
\begin{equation}
s_{(1^r)}^\se =
\begin{cases}
  s_{(1^r)} & \text{for $\se=\nn,\hdom\,$} \\
  s_{(1^r)}-s_{(1^{r-2})} & \text{for $\se=\vdom\,$} \\
  s_{(1^r)}-s_{(1^{r-1})} & \text{for $\se=\cell\,$.}
\end{cases}
\end{equation}

\subsection{Kernels} The Cauchy element $\Omega = \sum_{n \geq 0} s_n$
naturally plays a role in the generating function for the
structure coefficients of the ring of symmetric functions.
Since $s_\la[X+Y] = \sum_{\mu,\nu} c_{\mu\nu}^\la s_\mu[X]
s_\nu[Y]$ and
\begin{equation}\label{reproducing}
\Omega[XY] = \sum_{\la} s_\la[X] s_\la[Y],
\end{equation}
 then putting these two
identities together we have a generating function for all of the
LR coefficients:
\begin{equation*}
\Omega[X(Y + Z)] = \sum_{\la,\mu,\nu} c^{\la}_{\mu\nu} s_{\la}[X]
s_\mu[Y] s_\nu[Z]
\end{equation*}

Similarly we have $s_\la[XY] = \sum_{\mu,\nu} k_{\la\mu\nu}
s_\mu[X] s_\nu[Y]$ where $k_{\la\mu\nu}$ are the Kronecker
coefficients arising the expression $s_\la \ast s_\mu = \sum_{\nu}
k_{\la\mu\nu} s_\nu$ and $\frac{p_\la}{z_\la} \ast
\frac{p_\mu}{z_\mu} = \delta_{\la\mu} \frac{p_\la}{z_\la}$.
Substituting $YZ$ for $Y$ in in equation (\ref{reproducing})
yields
\begin{equation*}
\Omega[XYZ] = \sum_{\la,\mu,\nu} k_{\la\mu\nu} s_{\la}[X] s_\mu[Y]
s_\nu[Z].
\end{equation*}

The structure coefficients $d_{\la\mu\nu}$ satisfy the following
property.

\begin{prop}
\begin{equation}
\Omega[XY + XZ+YZ] = \sum_{\la,\mu,\nu} d_{\la\mu\nu} s_\la[X]
s_\mu[Y] s_\nu[Z]
\end{equation}
\end{prop}
\begin{proof}
\begin{align*}
\Omega[XY + XZ+YZ] &= \sum_{\alpha, \beta, \tau} s_\alpha[X]
s_\alpha[Y] s_\beta[X]
s_\beta[Z] s_\tau[Y] s_\tau[Z]\\
&=\sum_{\alpha, \beta, \tau, \la,\mu,\nu} c_{\alpha\beta}^\la
c_{\alpha\tau}^\mu c_{\beta\tau}^\nu
s_\la[X] s_\mu[Y] s_\nu[Z]\\
&= \sum_{\la,\mu,\nu} d_{\la\mu\nu} s_\la[X] s_\mu[Y] s_\nu[Z]
\end{align*}
\end{proof}

Both the left and right hand side of the above expression
are independent of $\se \in \{ \hdom\,, \vdom\,, \cell \}$ and yet
the coefficients $d_{\la\mu\nu}$ appear in this expression.
It turns out that this is only because the kernel $\Omega[XY + XZ + YZ]$
is common to all three bases as we show in the following corollary.

\begin{lemma}
\begin{equation*}
f_\se[X+Y] = f_\se[X] + XY + f_\se[Y]
\end{equation*}
\end{lemma}
\begin{proof}
Notice that $s_{(1^2)}[X + Y] = s_{(1^2)}[X] + XY + s_{(1^2)}[Y]$, $s_2[X+Y] = s_2[X]
+ XY + s_2[Y]$ and $s_{1}[X + Y] + s_{(1^2)}[X+Y] = s_{1}[X] + s_{1}[Y]+ s_{(1^2)}[X]
+ XY + s_{(1^2)}[Y]$.  Therefore we have for all $\se \in \{
\hdom\,, \vdom\,, \cell \}$.
\end{proof}

\begin{cor} For $\se \in \{ \hdom\,, \vdom\,, \cell\}$,
\begin{equation*}
\Omega[f_\se[X+Y+Z]] = \sum_{\la,\mu,\nu} d_{\la\mu\nu}
s_\la^{\se\ast}[X] s_\mu^{\se\ast}[X] s_\nu^{\se\ast}[X]
\end{equation*}
\end{cor}
\begin{proof}
From the previous lemma we can derive that $\Omega[f_\se[X+Y]] =
\Omega[XY] \Omega[f_\se[X]] \Omega[f_\se[Y]] = \sum_{\la }
s_\la^{\se\ast}[X] s_\la^{\se\ast}[Y]$.
In addition we have,
\begin{align*}
\Omega[f_\se[X+Y+Z]] &= \Omega[f_\se[X] +f_\se[Y] + f_\se[Z]
+ XY+XZ+YZ]\\
&= \Omega[f_\se[X] +f_\se[Y] + f_\se[Z]]
\sum_{\la,\mu,\nu} d_{\la\mu\nu} s_\la[X] s_\mu[Y] s_\nu[Z]\\
&=\sum_{\la,\mu,\nu} d_{\la\mu\nu} s_\la^{\se\ast}[X]
s_\mu^{\se\ast}[Y] s_\nu^{\se\ast}[Z].
\end{align*}
\end{proof}

\section{Hall-Littlewood symmetric functions and analogues}
\label{sec:qanalogues}

\subsection{Deformed Schur basis}
Define
\begin{equation}
  \tB(Z) = \sum_{\nu\in\Z^n} z^\nu \tB_\nu
\end{equation}
where $\tB_\nu$ is the $t$-analogue of $B_\nu$ from equation
(\ref{tanal}). This is the ``parabolic modified" analogue of
Jing's Hall-Littlewood creation operator. It was studied in
\cite{SZ}, where it is denoted by $H^t_\nu$. By \eqref{eq:itBern}
and \eqref{eq:analogmultskew},
\begin{equation} \label{eq:tpar}
  \tB(Z) = R(Z) \Omega[ZX] \Omega[(t-1) Z^* X]^\perp.
\end{equation}
Let $Z^{(1)},\dotsc,Z^{(L)}$ be a family of finite ordered
alphabets and $R_1$ through $R_L$ partitions such that the number
of parts of $R_j$ is equal to the number of letters in $Z^{(j)}$
for all $j$. Define the symmetric functions $\tBB_R[X;t]$ and
polynomials $c_{\la;R}(t)$ by
\begin{equation} \label{eq:genKostkacoefs}
  \tB_{R_1} \dotsm \tB_{R_L} 1 = \tBB_R[X;t] = \sum_\la s_\la[X] c_{\la;R}(t).
\end{equation}
The $c_{\la;R}(t)$ are the generalized Kostka polynomials of
\cite{SW}, as proved in \cite{SZ}.

By \eqref{eq:analoguezero} and \eqref{eq:crSchur} we have
\begin{equation} \label{eq:BR0}
  \tBB_R[X;0] = B_{R_1} \dotsm B_{R_L} 1 = s_{(R_1,\dotsc,R_L)}[X]
\end{equation}
where $(R_1,\dotsc,R_L)$ denotes the sequence of integers obtained
by juxtaposing the parts of the partitions $R_j$. By
\eqref{eq:analogueone} and \eqref{eq:crSchur} we have
\begin{equation} \label{eq:BR1}
  \tBB_R[X;1] = s_{R_1}[X] \dotsm s_{R_L}[X].
\end{equation}

\begin{example} For $R = (\young{&&\cr}\,,
\young{&\cr&\cr}\,, \young{\cr}\,)$ we have
\begin{equation}
\begin{split}
\tBB_R[X;t] = &~\young{\cr&\cr&\cr&&\cr}~+~t~\young{&\cr&&\cr&&\cr}
~+~t~\young{\cr\cr&\cr&&&\cr}~+~(t^2+t)~\young{&\cr&\cr&&&\cr}~+~
t^2~\young{\cr&&\cr&&&\cr}\\
&+~(t^3+t^2)~\young{\cr&\cr&&&&\cr}~+~
t^3~\young{&&\cr&&&&\cr}~+~t^4~\young{&\cr&&&&&\cr}
\end{split}
\end{equation}
\end{example}

\subsection{Deformed $s^\se_\la$ basis} \label{sec:deformed}
Let $\tB^\se_\nu$ be the $t$-analogue of $B^\se_\nu$. For
$\se\in\{\cell\,,\vdom\,,\hdom\,\}$ define
\begin{equation}
  \tB^\se(Z) = \sum_{\nu\in\Z^n} z^\nu \tB^\se_\nu.
\end{equation}
By \eqref{eq:analogmultskew}, Proposition \ref{pp:dbern},
\eqref{eq:Bcell} and \eqref{eq:Bhdom},
\begin{equation} \label{eq:tB}
\begin{split}
\tB^\vdom(Z) &= R(Z) \Omega[-s_{(1^2)}[Z]] \Omega[ZX]
\Omega[(Z+Z^*)(t-1)X]^\perp \\
\tB^\cell(Z) &= \tB^\vdom(Z) \Omega[-Z] \\
\tB^\hdom(Z) &= \tB^\vdom(Z) \Omega[-(1+\minus)Z].
\end{split}
\end{equation}
For a sequence of partitions $R=(R_1,R_2,\dotsc,R_L)$, define the symmetric function
$\tBB^\se_R[X;t]$ and the polynomials $d^\se_{\la R}(t)$ by
\begin{equation}
\tBB^\se_R[X;t] = \tB^\se_{R_1} \tB^\se_{R_2} \dotsm
\tB^\se_{R_L} 1 = \sum_\la d^\se_{\la R}(t) s^\se_\la.
\end{equation}
By \eqref{eq:analoguezero}, \eqref{eq:analogueone}, \eqref{eq:BR0}, \eqref{eq:BR1}, and
\eqref{eq:dberndef} we have
\begin{align}
\label{eq:rowspec0}
  \tBB^\se_R[X;0] &= s^\se_{(R_1,R_2,\dotsc,R_L)}[X] \\
\label{eq:rowspec1}
  \tBB^\se_R[X;1] &= s^\se_{R_1}[X] s^\se_{R_2}[X] \dotsm
  s^\se_{R_L}[X].
\end{align}

\begin{theorem} $d^\se_{\la R}(t)$ is constant over
$\se\in \{ \cell\,,\vdom\,,\hdom\,\}$.
\end{theorem}
\begin{proof} This follows from the fact that the $t$-analogue operation
commutes with the change of basis operations between the kinds
$\{\cell\,,\vdom\,,\hdom\,\}$, when acting on $\tB^\vdom(Z)$. For
by Propositions \ref{pp:dchange} and \ref{pp:dbern} we have
\begin{align}
  \Omega[-X]^\perp \tB^\vdom(Z) \Omega[X]^\perp &=
  \Omega[-Z]\tB^\vdom(Z)=\tB^\cell(Z) \\
  \Omega[-p_2[X]]^\perp \tB^\vdom(Z) \Omega[p_2[X]]^\perp &=
  \Omega[-(1+\minus)Z] \tB^\vdom(Z) = \tB^\hdom(Z).
\end{align}
\end{proof}

Let us call these polynomials $d_{\la R}(t)$. An important special
case is when $R$ consists of single-rowed rectangles of sizes
given by the partition $\mu$; in this case write $d_{\la\mu}(t)$
instead of $d_{\la R}(t)$.

\begin{theorem} $d_{\la\mu}(t)\in\NN[t]$.
\end{theorem}

These polynomials appear again in Proposition \ref{pp:singlerow} and it
is a consequence of this proposition and equation (\ref{eq:skewrec}) that
provides a proof of this theorem.

\begin{example}  We present an example of the symmetric functions
created by these operators. Let $\mu=(3,2,1)$.  For $\se \in \{
\cell\,, \vdom\,, \hdom\, \}$, we will represent the function
$s_\la^\se$ by the diagram for the partition $\la$ superscripted
by $\se$.  Notice that our example is independent of $\se$.
\begin{equation*}
\begin{split}
\tBB_\mu^\se[X;t] ~=~
&\young{\cr&\cr&&\cr}^\se+t\young{&&\cr&&\cr}^\se+t \young{\cr\cr&&&\cr}^\se+(t^2+t) \young{&\cr&&&\cr}^\se+(t^2+t^3)\young{\cr&&&&\cr}^\se\\
&+t^4 \young{&&&&&\cr}^\se+(t^2+t) \young{&\cr&\cr}^\se +t\young{\cr\cr&\cr}^\se
+(2 t^2+t+t^3) \young{\cr&&\cr}^\se \\
&+(t^4+t^2+t^3) \young{&&&\cr}^\se+(t^2+t^3) \young{\cr\cr}^\se +(t^4+t^2+t^3) \young{&\cr}^\se
+t^4
\end{split}
\end{equation*}
One might expect that $d_{\la R}(t)\in\NN[t]$ if the sequence
$(R_1,R_2,\dotsc,R_L)$ is a partition, since this conjecture has
been made for $\se = \nn$ \cite{SZ}. An example where this fails
to be the case is $R_1=(3)$, $R_2=(2,2)$ and $R_3 = (1)$ and $\se
= \hdom\,$. There the coefficient of $s_{(1,1)}^\hdom$ is
$t^5+t^3-t^4$.
%Is is possible that the generalized conjecture
%is true if $(\ell(\nu^{(1)}),
%\ell(\nu^{(2)}), \ldots, \ell(\nu^{(k)}))$ is a partition
%and $(\nu^{(1)},\nu^{(2)},\ldots,\nu^{(k)})$ is a partition?
\end{example}

We now give a more explicit formula for $d_{\la\mu}(t)$.
Iterating in a manner similar to \eqref{eq:twobern}, we have
\begin{equation}\label{eq:genKostka}
  \tB^\vdom(z_1) \dotsm \tB^\vdom(z_n) =
  \tB^\vdom(z_1,\dotsc,z_n) \prod_{1\le i<j\le n}
  \Omega[t z_j(z_i+z_i^*)]
\end{equation}
Let us apply both sides to the element $1\in\La$ and take the
coefficient at $Z^\mu$. Let $n$ be the length of
the partition $\mu$ and let
$\Pi$ be the above product over $i<j$. We have
\begin{equation*}
\begin{split}
\tBB^\vdom_\mu[X;t] &= R(Z) \Omega[-s_{(1^2)}[Z]] \Omega[ZX] \Pi
\coeff_{z^\mu} \\ &= R(Z) \Omega[-s_{(1^2)}[X]]^\perp \Omega[ZX] \Pi
\coeff_{z^\mu} \\ &= R(Z) \Pi \sum_{\la\in\Par^n} s_\la[Z]
s^\vdom_\la[X] \coeff_{z^\mu}.
\end{split}
\end{equation*}
Taking the coefficient of $s^\vdom_\la[X]$ for $\la\in \Par^n$ and using
$\rho$ and $J$ of type $A_{n-1}$ we have
\begin{equation*}
\begin{split}
d_{\la\mu}(t) &= z^{-\rho} J(z^\rho) s_\la[Z] \Pi \coeff_{z^\mu} \\
&= J(z^{\la+\rho}) \Pi \coeff_{z^{\mu+\rho}} \\
&= \left(\sum_{w\in S_n} \minus(w) z^{w(\la+\rho)} \Pi \right)\coeff_{z^{\mu+\rho}} \\
&= \sum_{w\in S_n} \left(\minus(w) \Pi \coeff_{z^{\mu+\rho-w(\la+\rho)}}\right).
\end{split}
\end{equation*}
Replacing each letter in $Z$ by its reciprocal, we have
\begin{equation} \label{eq:dqKostant}
  d_{\la\mu}(t) = \sum_{w\in S_n} \minus(w) \prod_{1\le i<j\le n}
  \Omega[t z_j^*(z_i+z_i^*)]
\coeff_{z^{w(\la+\rho)-(\mu+\rho)}}
\end{equation}
The product of geometric series can be expanded to obtain
a formula resembling Lusztig's $t$-analogue of weight multiplicity.
We also note that
\begin{equation*}
\begin{split}
d_{\la\mu}(t) &= \sum_{w\in S_n} \minus(w) z^{\mu+\rho} \prod_{1\le
i<j\le n}
  \Omega[t z_j^*(z_i+z_i^*)] \coeff_{z^{w(\la+\rho)}} \\
  &= J(z^{\mu+\rho} \prod_{1\le i<j\le n}
  \Omega[t z_j^*(z_i+z_i^*)])\coeff_{z^{\la+\rho}} \\
  &= J^{-1}(z^\rho) J\left(z^{\mu+\rho} \prod_{1\le i<j\le n}
  \Omega[t z_j^*(z_i+z_i^*)]\right) \coeff_{s_\la[Z]}.
\end{split}
\end{equation*}

\begin{comment}
Until now $d_{\la\mu}(t)$ has only been defined for $\la\in\Par^n$.
We may consistently define $d_{\la\mu}(t)$ for $\la,\mu\in P^+_{A_{n-1}}=\Z^n_{\ge}$ by
\begin{equation} \label{eq:drational}
  \sum_{\la\in P^+} d_{\la\mu}(t) s_\la[Z] = \pi\left( z^\mu \prod_{1\le i<j\le n}
  \Omega[t z_j^*(z_i+z_i^*)]\right)
\end{equation}
where $\pi$ is the isobaric divided difference (Demazure)
operator for the longest permutation $w_0$ of $S_n$,
acting on the $Z$ variables.

Note also that adding $(r^n)$ to both $\la$ and $\mu$ does not change
$d_{\la\mu}(t)$, so one need only compute $d_{\la\mu}(t)$ for $\la,\mu\in\Par^n$.

Let $\la^*=-w_0\la$ for $\la\in P$. Note that $s_\la[Z^*]=s_{\la^*}[Z]$
for all $\la\in P$. By linearity one has that
$\pi f(Z) = \sum_{\la\in P^+} f_\la s_\la[Z]$ if and only if
$\pi f(w_0(Z^*)) = \sum_{\la\in P^+} f_\la s_\la[Z^*]=\sum_{\la\in P^+} f_{\la^*} s_\la[Z]$.
Therefore
\begin{equation*}
  \sum_{\la\in P^+} d_{\la^*\mu^*}(t) s_\la[Z] =
  \pi \left( z^\mu \prod_{1\le i<j\le n} \Omega[t z_i(z_j+z_j^*)] \right) .
\end{equation*}
It follows that certain type $D_n$ affine Kazhdan-Lusztig polynomials
can be expressed in terms of the polynomials $d_{\la\mu}(t)$.
\end{comment}

\section{Parabolic Hall-Littlewood operators and analogues}
\label{sec:Hse} For each $\se\in\{\nn,\cell\,,\vdom\,,\hdom\,\}$
we define a variant of the same parabolic Hall-Littlewood creation
operator. These will be the creation operators for the universal
affine characters.

\subsection{Definition of operators}
Write $\tB^\se_{t^2}(Z)$ for $\tB^\se(Z)$ with $t$ replaced by
$t^2$. Let
\begin{equation}
\label{eq:crop} H^\se(Z) =\sum_{\nu\in\Z^k} z^\nu H^\se_\nu =
\Omega[\KK_\se[tX]-\KK_\se[X]]^\perp \tB_{t^2}(Z)
\Omega[\KK_\se[X]-\KK_\se[t X]]^\perp.
\end{equation}

\begin{prop} \label{pp:H} For
$\se\in\{\nn,\cell\,,\vdom\,,\hdom\,\}$,
\begin{equation}
  H^\se(Z) = \Omega[\KK_\se[t Z]] \tB^\se_{t^2}(Z).
\end{equation}
\end{prop}
\begin{proof} There is nothing to prove for $\se=\nn$.
Consider $\se=\vdom\,$. By Proposition \ref{pp:e2perpcomm} with
$W=t^2-1$,
\begin{equation*}
\begin{split}
&\quad \,\, \Omega[(t^2-1)s_{(1^2)}[X]]^\perp \Omega[ZX]
\Omega[(1-t^2)s_{(1^2)}[X]]^\perp \\
&= \Omega[ZX]\Omega[(t^2-1)s_{(1^2)}[Z]] \Omega[(t^2-1)ZX]^\perp.
\end{split}
\end{equation*}
Using \eqref{eq:tpar} and \eqref{eq:tB} we have
\begin{equation*}
\begin{split}
  H^\vdom(Z) &= \Omega[(t^2-1)s_{(1^2)}[X]]^\perp R(Z) \Omega[ZX]
  \Omega[(t^2-1)Z^* X]^\perp \Omega[(1-t^2)s_{(1^2)}[X]]^\perp \\
  &= R(Z) \Omega[ZX] \Omega[(t^2-1)s_{(1^2)}[Z]]
  \Omega[(t^2-1)(Z+Z^*)X]^\perp \\
  &=\Omega[t^2 s_{(1^2)}[Z]] \tB^\vdom_{t^2}(Z).
\end{split}
\end{equation*}
Next let $\se=\hdom\,$. We have
\begin{equation*}
\begin{split}
  H^\hdom(Z) &= \Omega[(t^2-1)p_2[X]]^\perp H^\vdom(Z)
  \Omega[(1-t^2)p_2[X]]^\perp \\
  &= \Omega[(t-1)(1+\minus)X]^\perp \Omega[(t^2-1)s_{(1^2)}[Z]] R(Z) \Omega[ZX] \\
  & \hphantom{xxxxx} \Omega[(t^2-1)(Z+Z^*)X]^\perp \Omega[(1-t)(1+\minus)X]^\perp \\
  &= \Omega[(t^2-1)s_{(1^2)}[Z]] R(Z) \\
  &\hphantom{xxxxx}\Omega[Z(X+(t-1)(1+\minus))] \Omega[(t^2-1)(Z+Z^*)X]^\perp \\
  &= \Omega[t^2s_2[Z]] \Omega[-p_2[Z]] \tB^\vdom_{t^2}(Z) \\
  &= \Omega[t^2s_2[Z]] \tB^\hdom_{t^2}(Z).
\end{split}
\end{equation*}
Finally, for $\se=\cell\,$ we have
\begin{equation*}
\begin{split}
  H^\cell(Z) &=
  \Omega[(t-1)X]^\perp H^\vdom(Z) \Omega[(1-t)X]^\perp \\
  &= R(Z) \Omega[t^2 s_{(1^2)}[Z]] \Omega[Z(X+(t-1))] \Omega[(t^2-1)(Z+Z^*)X]^\perp \\
  &= \Omega[(t-1)Z] \Omega[t^2 s_{(1^2)}[Z]] \tB^\vdom_{t^2}(Z) \\
  &= \Omega[tZ+t^2s_{(1^2)}[Z]] \tB^\cell_{t^2}(Z).
\end{split}
\end{equation*}
\end{proof}

Let $R=(R_1,R_2,\dotsc,R_L)$ be a sequence of partitions. For
$\se\in\{\nn,\,\cell\,,\,\vdom\,,\,\hdom\,\}$ define
$\HH^\se_R[X;t]$ and $K^\se_{\la;R}(t)$ by
\begin{equation}
\label{eq:qch} \HH^\se_R[X;t]=\sum_\la K^\se_{\la;R}(t)\,
s^\se_\la[X] = H^\se_{R_1} H^\se_{R_2} \dotsm H^\se_{R_L} 1.
\end{equation}
{}From the corresponding properties of the operator $\tB(Z)$, one
obtains the specializations at $t=0$ and $t=1$, for all $\se$.
\begin{align}
\label{eq:at0} \HH^\se_R[X;0] &= s^\se_{(R_1,R_2,\dotsc,R_L)}[X] \\
\label{eq:at1} \HH^\se_R[X;1] &= s_{R_1}[X] s_{R_2}[X] \dotsm
s_{R_L}[X].
\end{align}

\begin{remark} \label{rem:prodSchur}
For any $\se$, by \eqref{eq:at1} $\HH^\se_R[X;t]$ is a
$t$-deformation of the product of Schur functions. Note also that
$K^\nn_{\la;R}(t) = c_{\la;R}(t^2)$ where the latter is a
generalized Kostka polynomial; see \eqref{eq:genKostkacoefs}.
\end{remark}

\subsection{$\HH^\se$ in terms of $\HH^\nn$}
Let $|R|=\sum_i |R_i|$. Observe that
\begin{equation*}
\HH^\se_R[X;t] = \Omega[\KK_\se[t X]-\KK_\se[X]]^\perp
\HH^\nn_R[X;t].
\end{equation*}
{}From this and the fact that $K^\nn_{\tau;R}(t)=0$ unless
$|R|=|\tau|$, it follows that for
$\se\in\{\nn,\cell\,,\vdom\,,\hdom\,\}$,
\begin{equation} \label{eq:skewrec}
  K^\se_{\la;R}(t) = t^{|R|-|\la|} \sum_{\substack{\tau\in\Par \\
  |\tau|=|R|}} K^\nn_{\tau;R}(t) \sum_{\substack{\mu\in\Par^\se \\ |\mu|=|R|-|\la|}}
  c^\tau_{\la\mu}.
\end{equation}

\begin{remark} \label{rem:onerect}
Let $R=(R_1)$ be a single rectangle. For $\mu\subset R_1$ let $\widetilde{\mu}$ be the partition obtained
by the 180 degree rotation of the skew shape $R_1/\mu$. Then $c^{R_1}_{\la\mu}=0$ unless $\la=\widetilde{\mu}$,
in which case the coefficient is $1$. We have
\begin{equation}
  \HH^\se_{(R_1)}[X;t] = \sum_{\substack{\mu\in\Par^\se \\
  \mu\subseteq R_1 }} t^{|\mu|}
  s^\se_{\widetilde{\mu}}[X].
\end{equation}
\end{remark}

The known transpose symmetry of the generalized Kostka polynomials
for sequences of rectangles, induces a symmetry for all the
polynomials $K^\se_{\la;R}(t)$. Let $||R|| = \sum_{i<j} |R_i \cap
R_j|$, $\nn\,^t=\nn$, $\cell\,^t=\cell\,$, $\vdom\,^t=\hdom\,$,
and $\hdom\,^t=\vdom\,$.

\begin{prop} \label{pp:transpose} Let $R$ be a dominant sequence
of rectangles (that is, one whose widths weakly decrease)
and $R'$ a dominant rearrangement of $R^t$. Then for
all partitions $\la$,
\begin{equation} \label{eq:transpose}
  K^{\se^t}_{\la^t;R'}(t) = t^{2(||R||+|R|-|\la|)}
  K^{\se}_{\la;R}(t^{-1}).
\end{equation}
\end{prop}
\begin{proof} For $\se=\nn$ the formula holds by combining
\cite{SZ}, which connects the creation operator formula with the
graded character definition of \cite{SW}, \cite{Sh}, which
connects the definition of \cite{SW} with a tableau formula, and
either \cite{Sh2} or \cite{ScWa}, which show that the tableau
formula satisfies the above transpose symmetry. In the other cases
the formula holds by the case $\se=\nn$, \eqref{eq:skewrec}, and
\eqref{eq:LRtr}.
\end{proof}

\subsection{Connection between $\tBB^\se$ and $\HH^\se$}
We now make explicit the connection between the polynomials
$d_{\la R}(t)$ and $K^\se_{\la;R}(t)$ that is implied by
Proposition \ref{pp:H}.
\begin{theorem} \label{th:HBconnect} For $\se \in \{\nn\,, \cell\,, \vdom\,, \hdom\}$,
and $R=(R_1, R_2, \ldots, R_k)$ a sequence of partitions,
\begin{equation}\label{eq:relHandB}
\HH^\se_R[X;t] = \sum_{\nu^{(i)} \in \Par^\se} \sum_{\gamma^{(i)}} t^{|\nu|}
\left( \prod_i c_{\nu^{(i)}\gamma^{(i)}-(a_i^{b_i})}^{R_i-(a_i^{b_i})} \right) \tBB^\se_\gamma[X;t^2]
\end{equation}
where the sum is over all sequences of partitions $\nu =
(\nu^{(1)}, \nu^{(2)}, \ldots, \nu^{(k)})$, $b_i = \ell(R_i)$,
and $\gamma = (\gamma^{(1)}, \gamma^{(2)}, \ldots, \gamma^{(k)})$ is
a sequence of $\gamma^{(i)} \in \Z^{b_i}$ with $\gamma^{(i)}$ weakly
decreasing and $a_i = \gamma^{(i)}_{b_i}$.
\end{theorem}
\begin{proof}  Recall that we have the relation $\tB_{\nu}^\se =
\sign(w) \tB_{w(\nu +\rho)-\rho}^\se$ for $w \in S_n$ and for every
$\nu \in \Z^n$ there is a $w \in S_n$ such that $w(\nu + \rho)-\rho = \alpha \in
\Z^n$
with $\alpha_1 \geq \alpha_2 \geq \cdots \geq \alpha_n$. Therefore
\begin{equation}
\tB^\se(Z) = \sum_\nu z^\nu \tB_\nu^\se = \sum_{\alpha} \sum_{w \in S_n}
\sign(w) z^{w(\alpha+\rho)-\rho} \tB_\alpha^\se
\end{equation}
where the sum on the right is over all $\alpha \in \Z^n$ with $\alpha$ a weakly decreasing
sequence.  From equation $R(Z) = z^{-\rho} J(z^\rho)$ and hence
\begin{equation*}
R(Z)^{-1} \tB^\se(Z) = \sum_{\alpha} J(z^\rho)^{-1} J(z^{\alpha+\rho}) \tB_\alpha^\se
\end{equation*}
Now if $\alpha_n = r$ then $z^{w(\alpha+\rho)-\rho} = z^{(r^n)}
z^{w(\alpha-(r^n)+\rho)-\rho}$ where $(r^n)$ is shorthand notation
for the sequence with $r$ repeated $n$ times. Therefore
\begin{equation*}
R(Z)^{-1} \tB^\se(Z) = \sum_{\alpha} z^{(\alpha_n^n)}
s_{\alpha-(\alpha_n^n)}[Z] \tB_\alpha^\se.
\end{equation*}
We know by proposition \ref{pp:H} that $H^\se(Z) = \Omega[f_\se[tZ]] \tB^{\se}_{t^2}(Z)$
and since $\Omega[f_\se[tZ]] = \sum_{\la \in \PPP^\se} t^{|\la|} s_\la[Z]$ we have
that
\begin{equation*}
R(Z)^{-1} H^\se(Z) = \sum_{\alpha} \sum_{\nu \in \PPP^\se} \sum_\gamma
z^{(\alpha_n^n)} t^{|\nu|} c_{\alpha-(\alpha_n^n),\nu}^{\gamma} s_{\gamma}[Z] \tB_\alpha^{t^2\se}.
\end{equation*}
The sum over $\gamma$ in this expression will be over all $\gamma \in \PPP$
with  the length of $\gamma$ less than or equal to $n$.  Multiplying both
sides of this equation by $R(Z)$ and using again the identity that
$R(Z) s_\la[Z] = \sum_{w \in S_n} \sign(w) z^{w(\la+\rho)-\rho}$ yields
\begin{equation*}
H^\se(Z) =\sum_{\alpha} \sum_{\nu \in \PPP^\se} \sum_\gamma \sum_{w \in S_n}
t^{|\nu|} c_{\alpha-(\alpha_n^n),\nu}^{\gamma} \sign(w)
z^{w(\gamma+\rho+(\alpha_n^n)) - \rho} \tB_\alpha^{t^2\se}.
\end{equation*}
Since the sum is over $\gamma \in \PPP$, only for $w$ equal to the identity
will $w(\gamma+\rho+(\alpha_n^n)) - \rho$ be a partition
and the coefficient of $z^R$ for $R \in \PPP$ will be
\begin{equation*}
H^\se_R = \sum_{\alpha} \sum_{\nu \in \PPP^\se} t^{|\nu|}
c_{\alpha-(\alpha_n^n),\nu}^{R-(\alpha_n^n)} \tB_\alpha^{t^2\se}
\end{equation*}
The theorem follows since for a sequence of rectangles
$R=(R^{(1)}, R^{(2)}, \cdots, R^{(k)})$ with $R^{(i)} \in \PPP$ we have
$\HH_R^\se[X;t] = H_{R^{(1)}}^\se H_{R^{(2)}}^\se \cdots H_{R^{(k)}}^\se 1$
and iterating the above formula yields (\ref{eq:relHandB}).
\end{proof}
\begin{prop} \label{pp:singlerow} Let $\mu$ be a partition and $R$ the sequence
of single-rowed partitions $(\mu_i)$. Then
\begin{equation} \label{eq:singlerow}
\HH^\vdom_R[X;t]=\tBB^\vdom_R[X;t^2]
\end{equation}
or equivalently,
\begin{equation}
  K^\vdom_{\la;R}(t) = d_{\la R}(t^2).
\end{equation}
\end{prop}
\begin{proof} This follows from Proposition \ref{pp:H} and
the fact that $s_{(1^2)}[z]=0$ for $z$ a single letter.
\end{proof}

\begin{example}
We choose $R = (\young{&\cr&\cr}\,, \young{\cr})$
as a sequence of rectangles that will index an example of the
$4$ types of symmetric functions that we define here.  The first
example $H^\nn_{R}[X;t]$ is a generating function for generalized
Kostka polynomials; see Remark \ref{rem:prodSchur}.

\begin{equation*}
\HH^\nn_{\left(\young{&\cr&\cr}\,, \young{\cr}\right)}[X;t] = \young{\cr&\cr&\cr} + t^2 \young{&\cr&&\cr}
\end{equation*}

%%%H_{(2,2),(1)}^\vdom[X;t] =
%%%(t^2+1)*t^2*s[2,1]+t^2*s[1,1,1]+t^2*s[3,2]+s[2,2,1]+t^4*(t^2+1)*s[1]
\begin{align*}
\HH^\vdom_{\left(\young{&\cr&\cr}\,, \young{\cr}\right)}[X;t] =
\young{\cr&\cr&\cr}^\vdom+t^2 \young{&\cr&&\cr}^\vdom+
t^2 \young{\cr\cr\cr}^\vdom+(t^2+t^4) \young{\cr&\cr}^\vdom
+(t^4+t^6) \young{\cr}^\vdom
\end{align*}

%%%H_{(2,2),(1)}^\vdom[X;t] =
%%%t^2*s[3,2]+s[2,2,1]+t^4*s[3]+(t^2+1)*t^2*s[2,1]+t^4*(t^2+1)*s[1]
\begin{align*}
\HH^\hdom_{\left(\young{&\cr&\cr}\,, \young{\cr}\right)}[X;t] &=
\young{\cr&\cr&\cr}^\hdom+t^2 \young{&\cr&&\cr}^\hdom
+(t^2+t^4) \young{\cr&\cr}^\hdom
+t^4 \young{&&\cr}^\hdom+(t^4+t^6) \young{\cr}^\hdom
\end{align*}

%%%H_{(2,2),(1)}^\cell[X;t] =
%%%s[2,2,1]+t^2*s[3,2]+
%%%+t^3*s[3,1]+(t^2+1)*t*s[2,2]+t*s[2,1,1]+
%%%t^4*s[3]+2*(t^2+1)*t^2*s[2,1]+t^2*s[1,1,1]
%%%+t^3*(2*t^2+1)*s[2]+t^3*(t^2+2)*s[1,1]
%%%+2*t^4*(t^2+1)*s[1]+t^5*(t^2+1)

\begin{align*}
\HH^\cell_{\left(\young{&\cr&\cr}\,, \young{\cr}\right)}[X;t] &=
\young{\cr&\cr&\cr}^\cell+t^2 \young{&\cr&&\cr}^\cell
+t \young{\cr\cr&\cr}^\cell+
(t+t^3) \young{&\cr&\cr}^\cell + t^3 \young{\cr&&\cr}^\cell\\
&~~~+t^2 \young{\cr\cr\cr}^\cell+ 2 (t^2+t^4) \young{\cr&\cr}^\cell
 + t^4 \young{&&\cr}^\cell
+(t^3+2 t^5)\young{&\cr}^\cell\\
&~~~+(2 t^3 + t^5) \young{\cr\cr}^\cell
+2 (t^4+t^6) \young{\cr}^\cell + (t^5 + t^7)
\end{align*}

Next we list the functions $\tBB^\se_\gamma[X;t^2]$ that are relevant to
the computation of $\HH_R^\se[X;t]$ for $\se \in \{ \vdom\,, \hdom\,, \cell \}$.
For $\se = \nn$ we have that $\tBB^\nn_R[X;t^2] = \HH^\nn_R[X;t]$.

\begin{equation*}
\tBB^\se_{\left((2,2),(1)\right)}[X;t^2] = t^2s_{{3,2}}^\se+s_{{2,2,1}}^\se+t^2s_{{2,1}}^\se
\end{equation*}

\begin{equation*}
\tBB^\se_{\left((2,1), (1)\right)}[X;t^2] = t^2s_{{2}}^\se+t^2s_{{1,1}}^\se+t^2s_{{3,1}}^\se+t^2s_{{2,2}}^\se+s_{{2,1,1}}^\se
\end{equation*}

\begin{equation*}
\tBB^\se_{\left((2,0), (1)\right)}[X;t^2] =
t^2s_{{3}}^\se+t^2s_{{2,1}}^\se+t^2s_{{1}}^\se
\end{equation*}

\begin{equation*}
\tBB^\se_{\left((1,1), (1)\right)}[X;t^2] =
t^2s_{{2,1}}^\se+s_{{1,1,1}}^\se+t^2s_{{1}}^\se
\end{equation*}

\begin{equation*}
\tBB^\se_{\left((1,0), (1)\right)}[X;t^2] =
t^2+t^2s_{{2}}^\se+t^2s_{{1,1}}^\se
\end{equation*}

\begin{equation*}
\tBB^\se_{\left((0,0), (1)\right)}[X;t^2] =
t^2s_{{1}}^\se
\end{equation*}

\begin{equation*}
\tBB^\se_{\left((2,-1), (1)\right)}[X;t^2] =
(t^2-1)s_{{2}}^\se
\end{equation*}

\begin{equation*}
\tBB^\se_{\left((1,-1), (1)\right)}[X;t^2] =
(t^2-1)s_{{1}}
\end{equation*}

\begin{equation*}
\tBB^\se_{\left((0,-1), (1)\right)}[X;t^2] =
(t^2-1)
\end{equation*}

\begin{equation*}
\tBB^\se_{\left(\la, (0) \right)}[X;t^2] = s_{\la}^\se
\end{equation*}

Observe that Theorem \ref{th:HBconnect} says we will have the following
relationship.

\begin{equation*}
\HH^\vdom_{\left(\cclccl,\cell\right)}[X;t] = \tBB^\vdom_{\left((2,2),(1)\right)}
+ t^2 \tBB^\vdom_{\left((1,1),(1)\right)}
+ t^4 \tBB^\vdom_{\left((0,0),(1)\right)}
\end{equation*}

\begin{equation*}
\HH^\vdom_{\left(\cclccl,\cell\right)}[X;t] = \tBB^\hdom_{\left((2,2),(1)\right)}
+ t^2 \tBB^\hdom_{\left((2,0),(1)\right)}
+ t^4 \tBB^\hdom_{\left((0,0),(1)\right)}
\end{equation*}

\begin{align*}
\HH^\cell_{\left(\cclccl,\cell\right)}[X;t] &= \tBB^\cell_{\left((2,2),(1)\right)}
+ t \tBB^\cell_{\left((2,1),(1)\right)}
+ t^2 \tBB^\cell_{\left((1,1),(1)\right)}
+ t^2 \tBB^\cell_{\left((2,0),(1)\right)}\\
&+ t^3 \tBB^\cell_{\left((1,0),(1)\right)}
+ t^4 \tBB^\cell_{\left((0,0),(1)\right)}
+ t \tBB^\cell_{\left((2,2),(0)\right)}
+ t^2 \tBB^\cell_{\left((2,1),(0)\right)}\\
&+ t^3 \tBB^\cell_{\left((1,1),(0)\right)}
+ t^4 \tBB^\cell_{\left((2,0),(0)\right)}
+ t^4 \tBB^\cell_{\left((1,0),(0)\right)}
+ t^5 \tBB^\cell_{\left((0,0),(0)\right)}\\
&+ t^3 \tBB^\cell_{\left((2,-1),(1)\right)}
+ t^4 \tBB^\cell_{\left((1,-1),(1)\right)}
+ t^5 \tBB^\cell_{\left((0,-1),(1)\right)}
\end{align*}

\end{example}

\section{$X=M=K$}
\label{sec:XMK} Consider a nonexceptional affine Lie algebra
$\gggg$ of type $X^{(r)}_N$, with canonical simple Lie subalgebra
$\gb$ of rank $n$. Let $\la$ be a partition whose number of parts
is sufficiently smaller than $n$. There is a standard way to view
$\la$ as a dominant integral weight of $\gb$, namely, that $\la$
is identified with the weight $\sum_{i=1}^{n-1}
(\la_i-\la_{i+1})\omega_i$ where $\omega_i$ is the $i$-th
fundamental weight of $\gb$ (with $n$ large enough so that no spin
weights occur). Let $R$ be a sequence of rectangles
$R=(R_1,R_2,\dotsc)$ where $R_i$ has $r_i$ rows and $s_i$ columns
where $r_i$ is sufficiently smaller than $n$. Let $B^{r_i,s_i}$ be
the crystal graph of the KR module $W_s^{(r)}$ \cite{HKOTT}
\cite{HKOTY}. Let $\Xb_{R,\la}(t)$ be the classically restricted
one-dimensional sum for the crystal
\begin{equation} \label{eq:B}
B^R = \bigotimes_i B^{r_i,s_i}
\end{equation}
at the isotypic component of the irreducible $U_q(\gb)$-crystal of
highest weight $\la$, using the normalized coenergy function. Let
$\Mb_{R,\la}(t)$ be the corresponding fermionic formula, denoted
by $M(R,\la,t^{-1})$ in \cite{HKOTT} \cite{HKOTY}.
In those two papers it is conjectured
that $\Xb=\Mb$. By analyzing the fermionic formulae one can show
the following.

\begin{lemma} \label{lem:afflimit}
There is a well-defined limiting polynomial
\begin{equation}
  \lim_{n\rightarrow\infty} \Mb_{R,\la}(t)
\end{equation}
as the rank $n$ goes to infinity. It depends only on $R$, $\la$, and
the affine family of $X^{(r)}_N$. Moreover, there are only four
distinct families of such polynomials, which shall be named as
follows.
\begin{enumerate}
\item For $A^{(1)}_n$: $\Mb^\nn_{R,\la}(t)$.
\item For $B^{(1)}_n$, $D^{(1)}_n$, and $A^{(2)}_{2n-1}$:
$\Mb^{\vdom}_{R,\la}(t)$.
\item For $C^{(1)}_n$ and $A^{(2)\dagger}_{2n}$: $\Mb^{\hdom}_{R,\la}(t)$.
\item For $D^{(2)}_{n+1}$ and $A^{(2)}_{2n}$:
$\Mb^{\cell}_{R,\la}(t)$.
\end{enumerate}
\end{lemma}
Note that the families are grouped according to the decomposition
of $B^{r,s}$ as a $U_q(\gb)$-crystal, or equivalently, according
to the attachment of the $0$ node. See the appendices of
\cite{HKOTY} \cite{HKOTT}.

\begin{conjecture} \label{conj:create}
For $R$ a dominant sequence of rectangles and for all
$\se\in\{\nn,\,\cell\,,\,\vdom\,,\,\hdom\,\}$,
\begin{equation}
\label{eq:affch} K^\se_{\la;R}(t) =
\Mb^{\se^t}_{R^t,\la^t}(t^{2/\epsilon})
\end{equation}
where $\epsilon=1$ except for $\se=\cell\,$ in which case
$\epsilon=2$.
\end{conjecture}

This conjecture is surprising: it proposes a simple
relationship (coming from the deformations of the type A
Hall-Littlewood creation operators)
between the type A universal affine characters, and those of all other
nonexceptional affine types. At $t=1$ this was essentially known \cite{Kl}.
However the formulae for the powers of $t$ occurring in the affine characters
given either by the fermionic formula or the one-dimensional sum,
do not at all suggest such a simple relationship. Perhaps the
virtual crystal methods of \cite{OSS} can be used to prove
Conjecture \ref{conj:create}.

Equation \eqref{eq:affch} holds for $\se=\nn$. The proof connects
five formulae for the generalized Kostka polynomials.
\begin{enumerate}
\item The modified Hall-Littlewood creation operator formula, or equivalently, the $K$
polynomials \cite{SZ}.
\item Their definition (call it $\chi_q$) as isotypic components of graded Euler characteristic
characters of modules supported in the nullcone \cite{SW}.
\item A tableau formula (call it $LR$) using Littlewood-Richardson tableaux and a
generalized charge statistic \cite{ScWa} \cite{Sh}.
\item The affine crystal theoretic formula $X$ given by the
one-dimensional sum of type $A_n^{(1)}$ \cite{ScWa} \cite{Sh3}.
\item The fermionic formula $M$ \cite{KR}.
\end{enumerate}
We have $K=\chi_q$ \cite{SZ}, $\chi_q=LR$ \cite{Sh}, $LR=X$
\cite{ScWa} \cite{Sh3}, and $X=M$ \cite{KSS}.

Equation \eqref{eq:affch} also holds for a single rectangle in all
nonexceptional affine types, due to the agreement of
\cite[Appendix A]{HKOTT} and \cite[Appendix A]{HKOTY} with the
formula in Remark \ref{rem:onerect}.

Observe that by combining Conjecture \ref{conj:create} and
Proposition \ref{pp:transpose} one obtains the following conjecture.

\begin{conjecture} \label{conj:transpose}
\begin{equation}
\Mb^\se_{R;\la}(t) = t^{\epsilon(||R|| + |R| - |\la|)}
\Mb^{\se^t}_{R^t,\la^t}(t^{-1}).
\end{equation}
\end{conjecture}
This was proved in \cite{KS} via a direct bijection for $\se=\nn$.
This is striking as it relates the fermionic formulae of different
types. This kind of relation is not at all apparent from the
structure of the fermionic formulae.

\end{document}